\documentclass[]{amsart}

\makeatletter
\@namedef{subjclassname@2020}{%
	\textup{2020} Mathematics Subject Classification}
\makeatother

\RequirePackage{amsmath,amsthm,amssymb}

\RequirePackage{tikz}

\RequirePackage{hyperref}

\RequirePackage{graphicx}

\numberwithin{equation}{section}

\usepackage[left=3.7cm,top=2.5cm,right=3.7cm]{geometry}              
\title{Linear Dynamical Systems}

\author[Gilmore]{Clifford Gilmore}

\date{}

\address{%
School of Mathematical Sciences, 
University College Cork, Ireland.
}  

\email{clifford.gilmore@ucc.ie}

\thanks{Support from the Irish Research Council via a Government of Ireland Postdoctoral Fellowship is gratefully acknowledged.}



\usepackage[UKenglish]{babel}

\usepackage{mathrsfs}

\usepackage[shortlabels]{enumitem}

\usepackage{mathtools, bbm}

\usepackage{tikz-cd}

\usetikzlibrary[positioning,decorations.pathreplacing,arrows.meta,patterns]

\usetikzlibrary{decorations.pathmorphing}
\usetikzlibrary{decorations.markings}
\usetikzlibrary{bending}

\pgfdeclarepatternformonly{my crosshatch dots}{\pgfqpoint{-1pt}{-1pt}}{\pgfqpoint{5pt}{5pt}}{\pgfqpoint{9pt}{9pt}}%
{
	\pgfpathcircle{\pgfqpoint{1pt}{1pt}}{.8pt}
	\pgfpathcircle{\pgfqpoint{4pt}{4pt}}{.8pt}
	\pgfusepath{fill}
}

\usepackage[labelformat=simple]{subcaption}

\newtheorem{theorem}{Theorem}[section]

\newcommand{\Lsp}[1]{\mathscr{L}(#1)}  		
\newcommand{\Lx}{\Lsp{X}}  		 			
\newcommand{\spn}[1]{\mathrm{span}\!\left\lbrace #1 \right\rbrace}

\newcommand{\cspn}[1]{\overline{\mathrm{span}}\!\left\lbrace #1 \right\rbrace} 

\newcommand{\wcomp}{W_{\psi,\varphi}}

\newcommand{\N}{\mathbb{N}} 
\newcommand{\C}{\mathbb{C}} 
\newcommand{\R}{\mathbb{R}} 
\newcommand{\Z}{\mathbb{Z}}
\newcommand{\D}{\mathbb{D}}  
\newcommand{\T}{\mathbb{T}}  
\newcommand{\B}{\mathbb{B}}  
\newcommand{\1}{\mathbbm{1}}  	
\newcommand{\entireFns}{H(\C)}  

\DeclareMathOperator{\udens}{\overline{dens}} 
\DeclareMathOperator{\ldens}{\underline{dens}} 
 
\DeclareMathOperator{\ubd}{\overline{Bd}}  

\newcommand{\hc}[1]{HC\left( #1 \right)}  

\newcommand{\returnSet}[2]{\mathcal{N}_T(#1, \, #2)}  
\newcommand{\nTU}{\returnSet{x}{U}}

\newcommand{\abs}[1]{\left| #1 \right|}
\newcommand{\norm}[1]{{\left\| #1 \right\|}}

\theoremstyle{definition}

\newtheorem{qu}{Question} 

\newlength{\xRadiusSmall}
\newlength{\xRadius}
\newlength{\yRadiusSmall}
\newlength{\yRadius}
\newlength{\xAxisRadius}
\newlength{\yAxisRadius}
\begin{document}
	
	\begin{abstract}
		This expository survey is dedicated to recent developments in the area of linear dynamics.  Topics include frequent hypercyclicity, $\mathcal{U}$-frequent hypercyclicity, reiterative hypercyclicity, operators of C-type, Li-Yorke and distributional chaos, and hypercyclic algebras.
	\end{abstract}
	
	\maketitle
	
	
	\section{Introduction}

	Chaos theory has been described in lay terms as the `science of surprises' and in everyday usage chaos typically depicts something wild or a state of disorder.  While this is adequate in ordinary parlance, it is natural to ask: what mathematically precise definition captures the essential properties of a chaotic dynamical system?
	
	By a \emph{dynamical system} we mean a pair $(X,T)$, where $X$ is a metric space and $T$ is a continuous map acting on $X$.  The investigation of dynamical systems is primarily concerned with the long term evolution of iterates of the map $T$, where $n$-fold iteration is denoted by
	\begin{equation*}
		T^n = T \circ T \circ \cdots \circ T, \quad n \geq 0. 
	\end{equation*}
	Devaney~\cite{Dev89} suggested that $(X,T)$ is chaotic if it possesses the following three characteristics.
	\begin{enumerate}[itemsep=.5ex, label=\arabic*.] 
		\item It cannot be simplified.  \label{TopTrans}
		
		\item It has some regularity.  \label{densePeriodic}
		
		\item Long term prediction is difficult.  \label{SDIC}
	\end{enumerate}
	He proposed these characteristics are captured by the following three mathematical properties.
	
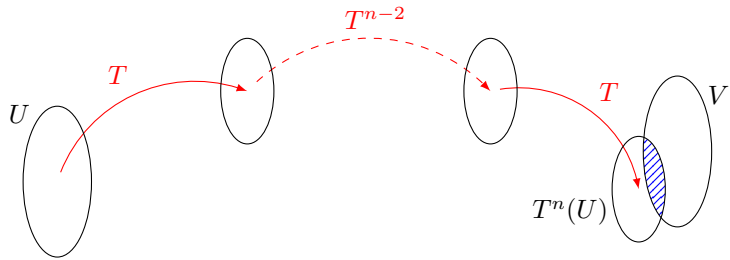
\begin{figure}[t]
		\centering
		\begin{tikzpicture}[scale=1.0]		
		\coordinate (Uellipse) at (0,0);
		\coordinate (U1ellipse) at (2.5,1.2);
		\coordinate (V1ellipse) at (5.7,1.2);
		\coordinate (Vellipse) at (8.16,0.4);
		\coordinate (TUellipse) at (7.65,-0.1);

		\node[] at (-0.5,0.9) (Ulabel) {$U$};
		\node[] at (8.7,1.15) (Vlabel) {$V$};
		\node[] at (6.75,-0.40) (Vlabel) {$T^n(U)$};

		\node[] at (Uellipse) (Unode) {}
		edge[-latex, bend left=45,color=red] node[auto] {$T$} (U1ellipse);
		
		\node[] at (U1ellipse) (U1node) {}
		edge[-latex,dashed, bend left=45,color=red] node[auto] {$T^{n-2}$} (V1ellipse);
		
		\node[] at (V1ellipse) (V1node) {}
		edge[-latex, bend left=45,color=red] node[auto] {$T$} (TUellipse);

		\begin{scope}
		\clip (Vellipse) ellipse [x radius=\xRadius,y radius=\yRadius,rotate=90];
		\fill[pattern=north east lines,pattern color=blue] (TUellipse) ellipse [x radius=\xRadiusSmall,y radius=\yRadiusSmall,rotate=90];
		\end{scope}
		
		\draw (Uellipse) ellipse [x radius=\xRadius,y radius=\yRadius,rotate=90];
		\draw (Vellipse) ellipse [x radius=\xRadius,y radius=\yRadius,rotate=90];
		\draw (TUellipse) ellipse [x radius=\xRadiusSmall,y radius=\yRadiusSmall,rotate=90];
		
		\draw (U1ellipse) ellipse [x radius=\xRadiusSmall,y radius=\yRadiusSmall,rotate=90];
		\draw (V1ellipse) ellipse [x radius=\xRadiusSmall,y radius=\yRadiusSmall,rotate=90];
		
		\end{tikzpicture}
		\caption{Topological transitivity.}
		\label{fig:topTrans}
	\end{figure}

	The first characteristic of chaos corresponds to the notion of topological transitivity.
	A dynamical system $(X,T)$ is said to be \emph{topologically transitive} if for any pair of nonempty, open subsets $U, V \subset X$, there exists some $n \in \mathbb{N}$ such that 
	\begin{equation*}
		T^n(U) \cap V \neq  \varnothing.
	\end{equation*}
	As illustrated in Figure \ref{fig:topTrans}, under the action of $T$ every non-trivial part of $X$ will eventually visit the whole space.   This captures how the system cannot be simplified or reduced into smaller and potentially more manageable components.  
	%
	%
	%
	%

	To satisfy regularity, Devaney defined the second characteristic to be when the map $T$ possesses a dense set of periodic points.  A vector $y \in X$ is a periodic point for $T$ if there exists $n \geq 1$ such that $T^n(y) = y$.  
	
	The third characteristic of chaos corresponds to the  notion of \emph{sensitive dependence on initial conditions}, which is commonly referred to as the \emph{butterfly effect}.  It is considered the essence of chaos since it describes how small discrepancies in the initial state of the system may lead to vastly different outcomes. According to E.~N.~Lorenz, who is considered the father of chaos theory, it characterises `when the present determines the future, but the approximate present does not approximately determine the future.'  It explains, for instance, the difficulty in obtaining accurate long-term weather forecasts.  
	
	An elegant result by Banks et al.~\cite{BBCDS92} demonstrated that sensitivity to initial conditions is redundant in Devaney's definition of chaos, since it follows automatically from the other two properties.  So we define a dynamical system $(X,T)$ to be \emph{chaotic} if $T$ is topologically transitive and $T$ possesses a dense set of periodic points.

	Chaos is typically  viewed as a nonlinear phenomenon.  Indeed, the field of classical dynamical systems investigates the mathematical rules governing the long-term evolution of nonlinear phenomena such as the weather, climate, turbulence and fluid dynamics.  
	However, it is now well established that seemingly tractable linear systems may give rise to complex dynamical behaviour and linear chaos.
	
	We say $(X,T)$ is a  \emph{linear dynamical system} if $X$ is a topological vector space and $T \colon X \to X$ is a continuous linear map.
	The central notion of linear dynamics is \emph{hypercyclicity}, since when $X$ is a separable, complete and metrizable topological vector space, the Birkhoff transitivity theorem gives that $T$ is hypercyclic if and only if $T$ is topologically transitive.

	Linear dynamics began to receive systematic attention in the early 1990s following the pioneering work of  Kitai~\cite{Kit82}, Gethner and Shapiro~\cite{GS87}, and Godefroy and Shapiro~\cite{GS91}.  It has since developed into a substantial branch of operator theory, which is evident from the monographs by Bayart and Matheron~\cite{BM09}, and Grosse-Erdmann and Peris~\cite{GEP11} that provide accessible and comprehensive introductions to the area.  
	
	This survey will primarily focus on advances that have occurred since publication of the books \cite{BM09} and \cite{GEP11}.  However, since this article is intended to be accessible to a wide mathematical audience, to ensure readability we begin by recalling the pertinent foundational concepts of linear dynamics.  Thus, before we highlight recent progress we set the scene by mentioning some significant background (or dare we say \emph{classical!}) results.  Finally, it is not possible to attempt an exhaustive account of the field and the topics selected here are entirely influenced by the personal preferences of the author.

	\section{Hypercyclicity}  \label{sec:Hc}
	
	For clarity of the presentation we will mostly consider the setting of Banach or Hilbert spaces.  Most of the definitions and results mentioned in the sequel can be appropriately extended to more general topological vector spaces such as Fr\'{e}chet spaces.  Where we consider it pertinent we will present results in this generality, so we remind the reader that a Fr\'{e}chet space is a locally convex and completely metrizable topological vector space.  
	So unless otherwise stated, we let $X$ be a separable Banach space and we denote by $\Lx$ the space of continuous linear operators on $X$.
	
	We say that $T \in \Lx$ is \emph{hypercyclic} if there exists $x \in X$ such that its $T$-orbit is dense in $X$, that is
	\begin{equation*}
		\overline{\{T^n x : n \geq 0\}} = X.
	\end{equation*}
	Such an $x \in X$ is called a \emph{hypercyclic vector} for $T$.
	
	As briefly mentioned in the introduction, the Birkhoff transitivity theorem~\cite{Bir20} states that an operator $T$ is hypercyclic if and only if it is topologically transitive.  
	This is a very useful result, and we note that in the course of the proof it  employs a Baire category argument to show that the set of hypercyclic vectors for a hypercyclic operator forms a dense $G_\delta$ subset of $X$  (cf.~\cite[Theorem 1.2]{BM09}).  We recall that a set is said to be $G_\delta$ if it is a countable intersection of open sets.
	
	A prerequisite for hypercyclicity is that the space $X$ is separable and furthermore hypercyclicity is a purely infinite-dimensional phenomenon, since linear operators cannot have dense orbits in the finite-dimensional setting. So we make a standing assumption that the spaces considered here are separable and infinite-dimensional.
	
	The first examples of hypercyclic operators in the Banach and Hilbert space settings were identified by  Rolewicz in 1969.  
	For the convenience of the reader we recall, for $1 \leq p < \infty$, that the space $\ell^p$ of $p$-summable sequences is a separable Banach space when endowed with the norm $\norm{x}_p = \left( \sum_{n=0}^\infty \abs{x_n}^p \right)^{1/p}$, for the sequence $x = (x_n)$. The space $c_0$ is defined as the space of sequences with limit equal to zero, which is a separable Banach space when endowed with the sup-norm $\norm{x} = \sup_{n \in \N} \abs{x_n}$, for  $x = (x_n)$. 
	Rolewicz~\cite{Rol69} proved in the setting $X = c_0$ or $\ell^p$, $1 \leq p < \infty$, that scalar multiples of  the backward shift $cB \in \Lx$ are hypercyclic when $\lvert c \rvert > 1$. The \emph{backward shift}
	$B \in \Lx$ is defined as
	\begin{equation*}
		B(x_1,\, x_2, \,x_3, \dotsc) = (x_2,\, x_3,\, x_4, \dotsc)
	\end{equation*}
	for $(x_n) \in X$.
	More generally, the \emph{weighted backward shift} $B_w \in \Lx$ is defined as 
	\begin{equation*}
		B_w(x_1,\, x_2, \, x_3, \dotsc ) = (w_2 x_2,\, w_3 x_3, \, w_4 x_4, \dotsc )
	\end{equation*}
	where $w = (w_n)$ is a bounded sequence of nonzero scalars.  For brevity we will simply refer to $B_w$ as a weighted shift.
	
	A complete characterisation of the hypercyclic weighted shifts was identified by Salas~\cite{Sal95}, who proved that $B_w$ is hypercyclic on $X$ if and only if 
	\begin{equation}  \label{charBwHc}
		\sup_{n\geq 1} \, \prod_{j=1}^n \abs{w_j} = \infty.
	\end{equation}
	In \cite{Sal95} it was also shown that any perturbation $I + B_w$ of the identity operator $I$ by a weighted shift $B_w$ is hypercyclic on $c_0$ or $\ell^p$, $1 \leq p < \infty$.
	
	The notion of hypercyclicity also makes sense in the setting of more general topological vector spaces.  In fact, the space $\entireFns$ of entire functions that are holomorphic on the complex plane provided the first examples of functions that admit dense orbits under linear operators.  We recall that $\entireFns$ is a Fr\'{e}chet space when endowed with the topology of local uniform convergence.  Birkhoff~\cite{Bir29} demonstrated in 1929 the existence of an entire function $f \in \entireFns$ such that the sequence of translates $\left(f(\, \cdot + na) \right)_{n\geq 1}$, for $a \neq 0$, forms a dense set in $\entireFns$.   MacLane~\cite{Mac52} subsequently constructed in 1952 an entire function $f \in \entireFns$ such that the sequence of derivatives $(f, \, f', \, f'', \dotsc)$ is dense in $\entireFns$.  We remark that the results from \cite{Bir29} and \cite{Mac52} were proven for the more general property of universality, which is briefly discussed in Section \ref{sec:motivation}.  In the language of linear dynamics these results give that translation operators $T_a \colon f(z) \mapsto f(z + a)$, for $a \neq 0$, and the differentiation operator $D \colon f \mapsto f'$ are hypercyclic on the space $\entireFns$.

	Some prevalent classes of operators that do not contain hypercyclic operators include contractions, finite-rank,  nuclear, compact, quasinilpotent, strictly singular and Riesz operators.  In the Hilbert space setting operators that are never hypercyclic include unitary, self-adjoint, normal, hyponormal, trace-class and Hilbert-Schmidt operators.  An account of families of non-hypercyclic operators can be found in \cite[Chapter 5]{GEP11}.
	
	Hypercyclicity is not, however, an anomalous phenomenon.  It was independently demonstrated by Ansari~\cite{Ans97} and Bernal~\cite{Ber99} that every separable, infinite-dimensional Banach space supports a hypercyclic operator.  This was generalised to the Fr\'{e}chet space setting by Bonet and Peris~\cite{BP98}, and Grivaux~\cite{Gri05} subsequently proved that the result holds for the stronger property of mixing.
	
	We say $T \in \Lx$ is \emph{mixing} if for any pair $U,\, V$ of nonempty open subsets of $X$, there exists $N \in \N$ such that $T^n U \cap V \neq \varnothing$ for all $n \geq N$.  Mixing is a strengthening of topological transitivity, so if $T$ is mixing the Birkhoff transitivity theorem gives that $T$ is hypercyclic.
	
	We gather the existence results into the following theorem.
	\begin{theorem}[Ansari~\cite{Ans97}, Bernal~\cite{Ber99}, Bonet and Peris~\cite{BP98}, Grivaux~\cite{Gri05}]
		Every separable, infinite-dimensional Fr\'{e}chet space supports a mixing, and hence hypercyclic, operator.
	\end{theorem}

	Since every separable, infinite-dimensional Banach space supports a hypercyclic operator, it is necessary that there exist hypercyclic operators that are compact perturbations of the identity, i.e.\ of the form
	\begin{equation*}
		I + K,
	\end{equation*}
	where $K$ is a compact (or even nuclear) operator.
	Such examples follow from results in \cite{Sal95}, and they are quite remarkable since individually neither compact operators nor the identity can be hypercyclic.  Another curious family of hypercyclic operators are the rank one perturbations of unitary operators acting on the Hilbert space $\ell^2$ that were constructed by Grivaux~\cite{Gri12}. This was a strengthening of a result by Shkarin~\cite{Shk10c}, who constructed a rank two perturbation of a unitary operator acting on a Hilbert space.

	Important classes of maps that contain hypercyclic operators include composition operators acting on function spaces such as the Hardy, Bergman and Dirichlet spaces.  Composition and weighted composition operators are defined, respectively, as
	\begin{equation*}
		f \mapsto C_\varphi f = f \circ \varphi, \qquad 	f \mapsto W_{\psi, \varphi} f = \psi \cdot f \circ \varphi
	\end{equation*}
	for fixed analytic maps $\varphi$ and $\psi$.  
	There exists a rich literature on the class of hypercyclic composition operators and an account of the fundamental results can be found in \cite[Chapter 4]{GEP11}.  To illustrate the theory we briefly mention the elegant characterisation of hypercyclic composition operators acting on the Hardy space $H^2(\D)$.  It was proven by Bourdon and Shapiro~\cite{BS90, BS97} that $C_\varphi$ is hypercyclic on $H^2(\D)$ if and only if $\varphi$ is an automorphism of the unit disc $\D$ with no fixed point in $\D$.  
	
	More recently Bayart~\cite{Bay10}, and Bayart and Charpentier~\cite{BC13} characterised the linear dynamical properties of $C_\varphi$ acting on the Hardy space $H^2(\B^d)$, where $\varphi$ a linear fractional map on the unit Euclidean ball $\B^d \subset \C^d$.
	Bonet and Doma\'{n}ski~\cite{BD12} gave the following characterisation of the hypercyclic composition operators acting on the space $\mathscr{A}(\Omega)$ of real analytic functions on an open subset $\Omega \subset \R^d$.  For $\varphi$ a real analytic self-map of $\Omega$, they proved that $C_\varphi$ is hypercyclic on $\mathscr{A}(\Omega)$ if and only if $\varphi$ is injective, $\varphi'$ is never singular and $\varphi$ is a runaway sequence on $\Omega$.  A map $\varphi \colon \Omega \to \Omega$ is said to be  \emph{runaway} if for every compact subset $K \subset \Omega$, there exists $n \in \N$ such that $\varphi^n(K) \cap K = \varnothing$.

	In contrast, it was shown in \cite{CG20} and \cite{FGJ18} that the Fock and Schwartz spaces do not even support, respectively, supercyclic weighted and unweighted composition operators.
	We say the operator  $T \in \Lx$ is \emph{supercyclic} if there exists $x \in X$ such that its projective $T$-orbit  is dense in $X$, that is
	\begin{equation*}
		\overline{ \left\lbrace \lambda T^n x : n \geq 0, \: \lambda \in \mathbb{C} \right\rbrace } = X.
	\end{equation*}

	Further significant classes of hypercyclic operators include adjoints of multipliers on spaces of holomorphic functions, generalisations of backward shifts acting on the Hardy and Bergman spaces (cf.~\cite[Chapter 4]{GEP11}), and Toeplitz operators acting on the Hardy space \cite{BL16}.
	In the setting of $\Lx$ and its separable ideals, hypercyclic properties of the class of elementary operators were investigated in \cite{Cha99} and \cite{BMP04}.  As noted in \cite{GST17,Gil19}, some interesting questions remain open regarding the hypercyclicity of commutator maps and generalised derivations acting on separable Banach ideals of $\Lx$.
	
	An introduction to linear dynamics would be incomplete without mentioning the Hypercyclicity Criterion. It is a sufficient condition for hypercyclicity that has emerged as a powerful tool, since for a given operator it is not always straightforward to explicitly identify a hypercyclic vector.  It was one of the principal results from Kitai~\cite{Kit82}, which was independently rediscovered by Gethner and Shapiro~\cite{GS87}.
	
	We say $T \in \Lx$  satisfies the Hypercyclicity Criterion if there exist dense subsets $X_0$, $Y_0 \subset X$, an increasing sequence $(n_k)$ of positive integers and maps $S_{n_k}\colon Y_0 \to X$, for $k \geq 1$, such that for any $x \in X_0$, $y \in Y_0$ the following hold as $k \to \infty$
	\begin{enumerate}[label=(\roman*),  itemsep=1ex]
		\item $T^{n_k} x \to 0$,
		\item $S_{n_k}(y) \to 0$, 
		\item $T^{n_k}S_{n_k}(y) \to y$. 
	\end{enumerate}
	Note that the maps of the sequence $S_{n_k}$ are not assumed to be self-maps of $Y_0$, linear or even continuous.
	If $T$ satisfies the Hypercyclicity Criterion then $T$ is hypercyclic (cf.~\cite[Theorem 1.6]{BM07}). An elegant argument to prove this result shows that under the assumptions of the Hypercyclicity Criterion, the operator $T$ is topologically transitive.

	If $T$ satisfies the Hypercyclicity Criterion for the full sequence $(n_k) = (k)$ of natural numbers, then $T$ is mixing (cf.~\cite[Remark 3.13]{GEP11}). 
	Moreover, B\`{e}s and Peris~\cite{BP99} proved that $T \in \Lx$ satisfies the Hypercyclicity Criterion if and only if $T$ is weakly mixing.
	We say $T \in \Lx$ is \emph{weakly mixing} if the direct sum $T \oplus T$ is hypercyclic on $X \oplus X$.   
	
	A long-standing open problem, originally posed in 1991 by Herrero~\cite{Her91},  asked if $T$ is hypercyclic does it follow that $T \oplus T$ is hypercyclic? (Equivalently, does $T$ satisfy the Hypercyclicity Criterion?)
	The question was resolved in the negative in 2006 by de la Rosa and Read~\cite{DLRR09}, who constructed a Banach space and a hypercyclic operator $T$ such that $T$ does not satisfy the (so-called!) Hypercyclicity Criterion.  Bayart and Matheron~\cite{BM07} subsequently identified a family of hypercyclic, non-weakly mixing operators in the setting of the classical Banach spaces $\ell^p$, for $1 \leq p < \infty$, which we note includes examples in the Hilbert space setting.  A simplification of the counterexample from \cite{BM07} can be found in \cite[Section 4.2]{BM09}.

	\subsection{Motivation}  \label{sec:motivation}
	
	While linear dynamics has developed into a substantial research area in its own right, we briefly mention here something about its origins.
	
	One motivation for investigating hypercyclic operators grew from the \emph{invariant subspace problem} and the study of cyclic operators.
	We say $T \in \Lx$ is \emph{cyclic} if there exists $x \in X$ (said to be a \emph{cyclic vector} for $T$) such that the closed linear span of its $T$-orbit is dense in $X$, that is
	\begin{equation*}
		\cspn{T^n x : n \geq 0 } = X.
	\end{equation*}
	The invariant subspace problem asks, given $T \in \Lx$, does there always exist a non-trivial, closed $T$-invariant subspace $W \subset X$?  The subspace $W$ is $T$-invariant if $T(W) \subset W$ and it is said to be non-trivial if $W \neq \{ 0 \}$ and $W \neq X$. 
	Clearly $T$ does not possess a non-trivial closed invariant subspace if and only if every nonzero $x \in X$ is a cyclic vector for $T$. 
	
	A counterexample to the invariant subspace problem was constructed in 1976 by Enflo~\cite{Enf87}.  Read~\cite{Rea85} subsequently identified an operator $T$, acting on the classical Banach space $\ell^1$, such that every nonzero  $x \in \ell^1$ is cyclic for $T$.  
	However, the invariant subspace problem remains an open question in the Hilbert space and reflexive Banach space settings.  We refer the interested reader to the monographs \cite{RR03} and \cite{CP11} for in-depth studies of this famous problem.
	
	The invariant subspace problem naturally led research activity to the analogous \emph{invariant subset problem} and the study of hypercyclic operators.  
	For $T \in \Lx$, the closure of the $T$-orbit of $x \in X$ is the smallest closed $T$-invariant subset that contains $x$.  Thus $T$ does not admit a non-trivial closed invariant subset if and only if every nonzero $x \in X$ is a hypercyclic vector for $T$.  Read~\cite{Rea88} identified such an example by constructing an operator $T \colon \ell^1 \to \ell^1$ such that every nonzero $x \in \ell^1$ is a hypercyclic vector for $T$.  By taking the closed linear span of the orbits, it follows that we have a counterexample to the invariant subspace problem.  For recent contributions to this topic we refer the curious reader to \cite{GR08} and \cite{GR14}.

	On the other hand, motivation to investigate hypercyclicity also stems from the more general notion of universality.  For topological spaces  $X$ and $Y$, the countable family $\left( T_n \right)_{n \in \mathbb{N}}$ of continuous maps $T_n \colon X \to Y$ is \emph{universal} if there exists $x \in X$ such that
	\begin{equation*}
		\overline{\left\lbrace T_n (x) : n \in \mathbb{N} \right\rbrace} = Y.
	\end{equation*}
	Such an  $x \in X$ is called a \emph{universal element} for $\left( T_n \right)$.
	If we let $X =Y$ be a topological vector space and we take the sequence $\left( T_n \right)$ to be the iterates of a single linear operator, then it follows that hypercyclicity is a particular instance of universality.
	
	The discovery of universal power series was credited in 1914 to Fekete~\cite{Pal14}.  He showed there exists a formal real power series $\sum_{j=1}^\infty a_j x^j$ on $[ -1, 1 ]$  with the property that for any continuous function $g \colon \left[ -1, 1 \right] \to \mathbb{R}$ with $g(0) = 0$, there exists an increasing sequence of positive integers $(n_k)$ such that 
	\begin{equation*}
		\sum_{j=1}^{n_k} a_j x^j \to g(x)
	\end{equation*}
	uniformly as $k \to \infty$.  The observation of Fekete can be further extended to a universal Taylor series  on all of $\mathbb{R}$ (cf.~\cite[Section 3a]{Gro99}).
	
	Many of the statements for hypercyclicity have analogues for the more general notion of universality.   In fact, the hypercyclicity results for the classical translation \cite{Bir29} and differentiation \cite{Mac52} operators acting on the space $\entireFns$ were originally proven for universality. However, some of the powerful tools used to investigate hypercyclic operators, for instance the  spectral techniques, are not available for universality.
	To learn more about this intriguing topic we refer the interested reader to the survey by Grosse-Erdmann~\cite{Gro99} and the article by Bayart et al.~\cite{BGENP08}.

	\subsection{Linear Chaos}
	
	Following the seminal paper of Godefroy and Shapiro~\cite{GS91}, the definition proposed by Devaney became the accepted definition of chaos in linear dynamics.  Alternative definitions of chaos also appear in the literature and we discuss some of them in Section \ref{sec:DistChaos}.
	
	We recall that the Birkhoff transitivity theorem gives that an operator is hypercyclic if and only if it is topologically transitive, which leads to the following definition.  We say that $T \in \Lx$ is \emph{chaotic} if the following hold:
	\begin{enumerate}[label=(\roman*),  itemsep=.5ex]
		\item $T$ is hypercyclic,
		
		\item $T$ possesses a set of periodic points that is dense in $X$.
	\end{enumerate}
	
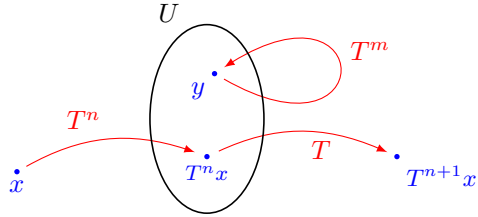
\begin{figure}[t]
		\centering
		\begin{tikzpicture}[scale=1]
		\coordinate (xVector) at (0,.5);
		\coordinate (U1ellipse) at (2.5,1.2);
		\coordinate (xEllipseVector) at (2.5,.7);
		\coordinate (awayVector) at (5,.7);
		\coordinate (x1) at (2.6, 0.7);
		\coordinate (y) at (2.6, 1.8);

		%
		\node at (y) (PeriodicPt) {};
		\fill [blue] (PeriodicPt) circle (1pt) node[below left] {\color{blue} $y$};
		\draw [-latex,color=red] (PeriodicPt) to [out=-30,in=34,looseness=40] node[above right] {$T^m$} (PeriodicPt);
		
		\node[] at (2.0,2.6) (Ulabel) {$U$};

		\node[] at (U1ellipse) (U1node) {};
		\node[] at (xEllipseVector) (xEllipseNode) {};
		\node[] at (awayVector) (awayNode) {};
		\node[] at (x1) (x1Node) {};

		\node[] at (xVector) (xNode) {}
		edge[-latex, shorten >=1pt, bend left=25,color=red] node[auto] {$T^n$} (xEllipseNode);

		\fill [blue] (xNode) circle (1pt) node[below] {\color{blue}$x$};

		\node[] at (awayVector) (awayNodeBelow) {}
		edge[latex-, shorten >=1pt, bend right=25,color=red] node[below right] {$T$} (xEllipseNode);

		\fill [blue] (xEllipseNode) circle (1pt) node[below] {\footnotesize{$T^n x$}};
		
		\fill [blue] (awayNode) circle (1pt) node[anchor=north west] {\color{blue}\small{$T^{n+1}x$}};
		
		\draw[line width=.6pt] (U1ellipse) ellipse [x radius=\xAxisRadius,y radius=\yAxisRadius,rotate=90];
		\end{tikzpicture}
		\caption{Chaos: each nonempty open subset contains a periodic point and a hypercyclic vector.}
		\label{fig:Chaos}
	\end{figure}

	The behaviour of a periodic point is in stark contrast to that of a hypercyclic vector.  However, as illustrated in Figure \ref{fig:Chaos}, for an operator $T$ to be chaotic each nonempty open subset $U \subset X$ must contain a hypercyclic vector and a periodic point.
	We also remark that for a linear dynamical system $(X,T)$, if $T$ is hypercyclic then it follows directly that $T$ has sensitive dependence on initial conditions \cite{GS91}. 
	

	The classical hypercyclic operators previously mentioned (non-trivial translations, differentiation operators and scalar multiples of the backward shift) all turn out to be chaotic (cf.~\cite[Section 2.3]{GEP11}).  
	Further examples of chaotic operators include adjoint multipliers and (weighted) composition operators acting on function spaces such as the Hardy and Bergman spaces \cite{GS91, Rez11, Bes14}. 
	In particular, the chaotic composition operators acting on the Hardy space $H^2(\D)$ were characterised by Taniguchi~\cite{Tan04}, who proved that $C_\varphi$ is chaotic on $H^2(\D)$ if and only if $\varphi$ is an automorphism of the unit disc $\D$ with no fixed point in $\D$.
	
	The chaotic weighted shifts $B_w$ acting on $\ell^p$, $1 \leq p < \infty$, were characterised by Grosse-Erdmann~\cite{Gro00}, who proved that $B_w$ is chaotic if and only if
	\begin{equation} \label{char:BwChaotic}
		\sum_{n=1}^\infty \frac{1}{\abs{w_1 \cdots w_n}^p} < \infty
	\end{equation}
	and if and only if $B_w$ admits a nonzero periodic point.  The weighted shift $B_w$ is chaotic on $c_0$ if and only if $\prod_{j=1}^n \abs{w_j} \to \infty$, as $n\to \infty$.

	A significant difference between hypercyclicity and chaos was identified by Bonet et al.~\cite{BMP01}, who demonstrated that there exist spaces that admit no chaotic operator.
	\begin{theorem}[Bonet, Mart{\'{\i}}nez-Gim{\'e}nez and Peris~\cite{BMP01}]  \label{thm:BMP01}
		There exist separable infinite-dimensional Banach spaces that do not support a chaotic operator.
	\end{theorem}
	To prove Theorem \ref{thm:BMP01}, they demonstrated that the hereditarily indecomposable Banach spaces constructed by Gowers and Maurey~\cite{GM93} do not support chaotic operators.  Thus operators of the form $I + K$ were shown to be non-chaotic, where $K$ is strictly singular and $I$ is the identity. 
	We note that their argument also holds for operators of the form $I + K$, where $K$ is compact, acting on the hereditarily indecomposable Banach spaces that were subsequently constructed by Argyros and Haydon~\cite{AH11}.  A Banach space $X$ is said to be hereditarily indecomposable if no closed subspace of $X$ is decomposable as a direct sum of infinite-dimensional subspaces.
	
	On the other hand, de la Rosa et al.~\cite{DLRLGP12} constructed chaotic operators that are compact (or even nuclear) perturbations of diagonal operators that have complex diagonal coefficients of modulus 1.  This gives, for instance, that there exist chaotic operators on any complex Banach space with an unconditional basis.
	
	\begin{theorem}[de la Rosa, Leonhard, Grivaux and Peris~\cite{DLRLGP12}]  \label{thm:RLGP12}
		Let $X$ be a complex separable Banach space having an unconditional Schauder decomposition. Then $X$ supports an operator which is chaotic.
	\end{theorem}
	
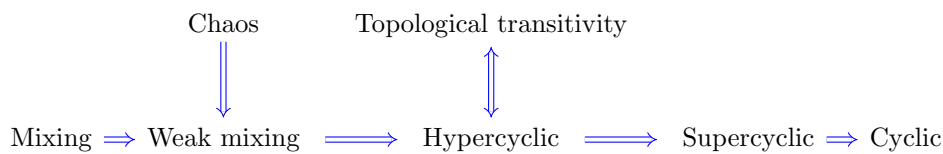
\begin{figure}[t]
		\centering
		\begin{tikzcd}[row sep=large, column sep=small, scale=1]
			& \textrm{Chaos} \arrow[d, Rightarrow, blue] &\textrm{Topological transitivity} \arrow[d, Leftrightarrow, blue] & &\\
			\textrm{Mixing} \arrow[r, Rightarrow, blue] & \textrm{Weak mixing} \ar[r,  Rightarrow, shorten >= .5em,shorten <= .5em, blue] & \textrm{Hypercyclic} \arrow[r, Rightarrow,shorten >= .5em,shorten <= .5em, blue] &  \textrm{Supercyclic} \arrow[r, Rightarrow, blue	] & \textrm{Cyclic}
		\end{tikzcd}	
		\caption{Relations between the dynamical properties from Section \ref{sec:Hc}.}
		\label{fig:sec2DynProps}
	\end{figure}

	Chaotic operators satisfy the Hypercyclicity Criterion and hence chaos implies weak mixing (cf.~\cite[Proposition 6.11]{BM07}).
	A summary of the relations between the dynamical properties introduced thus far can be found in Figure \ref{fig:sec2DynProps}.

	
	\section{Recurrence in Linear Dynamics}

	For $T \in \Lx$ to satisfy the definition of hypercyclicity, we require the existence of a vector $x \in X$ such that for every nonempty, open subset $U \subset X$, the \emph{return set} 
	\begin{equation*}
		\nTU \coloneqq \left\lbrace n \geq 0 \,:\, T^n x \in U \right\rbrace
	\end{equation*}
	is nonempty, as illustrated in Figure \ref{fig:HCdefn}.  It turns out if $T$ is hypercyclic then $\nTU$ is in fact an infinite set.
	It is thus natural to ask, are some orbits more recurrent than others?
	
	One approach to resolve this question is to consider quantitative differences in hypercyclic behaviour by calculating an appropriate density of the return sets.
	The \emph{lower} and \emph{upper densities} of a set $A \subset \N$ are defined, respectively, as
	\begin{align*}
		\ldens(A) &\coloneqq \liminf_{n \to \infty} \frac{ \abs{ A \cap \{ 1, 2, \dotsc , n\}}}{n}, \\
		\udens(A) &\coloneqq \limsup_{n \to \infty} \frac{ \abs{ A \cap \{ 1, 2, \dotsc , n\}}}{n},
	\end{align*}
	where $\abs{\,\cdot\,}$ denotes the cardinality of the set.
	The \emph{upper Banach density} of $A \subset \N$ is defined as
	\begin{equation*}
		\ubd(A) \coloneqq \lim_{N \to \infty} \frac{b_N}{N}
	\end{equation*}
	where
	\begin{equation*}
		b_N \coloneqq \limsup_{k \to \infty} \abs{ A \cap [k+1,\, k+N ]}.
	\end{equation*}
	These densities are related as follows
	\begin{equation*}
		\ldens(A) \leq \udens(A) \leq \ubd(A).
	\end{equation*}

	Investigation of the lower density of return sets was initiated in 2004 by Bayart and Grivaux~\cite{BG04,BG06}  when they introduced the notion of frequent hypercyclicity. Since then the study of recurrent orbits has developed into one of the most important branches of linear dynamics.

	\subsection{Frequent Hypercyclicity}
	
	We say $T \in \Lx$ is \emph{frequently hypercyclic} if there exists $x \in X$ such that for any nonempty open subset $U \subset X$, the return set $\nTU$ has positive lower density, that is
	\begin{equation*}
		\ldens(\nTU)  > 0.
	\end{equation*}
	Such an $x \in X$ is a \emph{frequently hypercyclic vector} for $T$.

	Frequent hypercyclicity gives a quantitative description of how frequently an orbit visits each neighbourhood.
	As illustrated in  \ref{fig:HC}, for $T$ to be frequently hypercyclic, each time the $T$-orbit visits a particular $U$ it must actually visit $U$ quite often before it moves away.


	\begin{figure}[t]
		\begin{subfigure}[]{.49\linewidth}
			\centering
			\begin{tikzpicture}[scale=1.1]
			\coordinate (xVector) at (0,1);
			\coordinate (U1ellipse) at (2.5,1.2);
			\coordinate (awayVector) at (5,1.2);
			\coordinate (x1) at (2.6, 0.7);
			\coordinate (x2) at (2.3, 0.4);
			\coordinate (x3) at (2.6, 2.0);
			\coordinate (x4) at (2.0, 1.2);
			\coordinate (x5) at (2.9, 1.3);

			\node[] at (3.1,2.6) (Ulabel) {$U$};

			\node[] at (xVector) (xNode) {}
			edge[-latex, shorten >=3pt, bend left=35,color=red] node[auto] {$T^n$} (U1ellipse);

			\fill [blue] (xNode) circle (1pt) node[below] {\color{blue}$x$};
			
			\node[] at (U1ellipse) (U1node) {};
			\node[] at (awayVector) (awayNode) {};
			\node[] at (x1) (x1Node) {};
			
			\node[] at (awayVector) (awayNodeBelow) {}
			edge[latex-, shorten >=3pt, bend right=35,color=red] node[above right] {$T$} (U1ellipse);
			
			\fill [blue] (U1node) circle (1pt) node[below] {\footnotesize{$T^n x$}};
			
			\fill [blue] (awayNode) circle (1pt) node[anchor=north west] {\color{blue}\small{$T^{n+1}x$}};
			
			\draw[line width=.6pt] (U1ellipse) ellipse [x radius=\xAxisRadius,y radius=\yAxisRadius,rotate=90];
			
			\end{tikzpicture}
			\caption{$T$-orbit of $x$ visits $U$.}
			\label{fig:HCdefn}
		\end{subfigure}		
		\hfill
		\begin{subfigure}{.49\linewidth}
			\centering
			\begin{tikzpicture}[scale=1.0]
			\coordinate (xVector) at (0,1);
			\coordinate (U1ellipse) at (2.5,1.2);
			\coordinate (x1) at (2.5, 0.9);

			\node[] at (3.5,0.6) (Ulabel) {$U$};

			\node[] at (4.2,3) (Ulabel) {\color{red}$T \circ T \circ \cdots$};

			\fill [blue] (xNode) circle (1pt) node[below] {\color{blue}$x$};

			\node[] at (U1ellipse) (U1node) {};
			\node[] at (awayVector) (awayNode) {};
			\node[] at (x1) (x1Node) {};
			
			\draw[line width=.6pt] (U1ellipse) ellipse [x radius=\xAxisRadius,y radius=\yAxisRadius,rotate=90];
			
			\fill[pattern=my crosshatch dots, pattern color=blue] (U1ellipse) ellipse [x radius=\xAxisRadius,y radius=\yAxisRadius,rotate=90];

			\node[] at (xVector) (xNode) {}
			edge[-latex, shorten >=3pt, bend left=25,color=red] node[above left] {$T^n$} (x1);

			\draw[decorate, decoration={coil,segment length = 9.9pt,aspect=0.3,amplitude=25pt},color=red,arrows = {latex-}] (2.5,3.0) -- (2.47,.9);
			\end{tikzpicture}
			\caption{Recurrent $T$-orbit of $x$ in $U$.}
			\label{fig:HC}
		\end{subfigure}
		\caption{}
		\label{fig:hypercyclicity}
	\end{figure}

	The classical hypercyclic operators mentioned in Section \ref{sec:Hc} (non-trivial translations, differentiation and scalar multiples of the backward shift) turn out to possess the stronger property of frequent hypercyclicity.  
	Weighted shifts also provide a rich source of examples to illustrate the nature of frequent hypercyclicity.

	Frequently hypercyclic weighted shifts were characterised by Bayart and Ruzsa~\cite{BR15} with the following theorem.
	
	\begin{theorem}[Bayart and Ruzsa~\cite{BR15}]  \label{thm:fhcBwlp}
		Let $1 \leq p < \infty$.  The weighted shift $B_w \colon \ell^p \to \ell^p$ is frequently hypercyclic if and only if
		\begin{equation} \label{charBwFhc}
			\sum_{n=1}^\infty \frac{1}{\abs{w_1 \cdots w_n}^p} < \infty.
		\end{equation}
	\end{theorem}
	We remark that a comparison with \eqref{char:BwChaotic} reveals that a weighted shift acting on $\ell^p$, for $1 \leq p < \infty$, is frequently hypercyclic if and only if it is chaotic.  
	
	Weighted shifts also provide examples of hypercyclic operators that are not frequently hypercyclic.
	For instance, in \cite[Example 6.17]{BM09} they show that $B_w \colon \ell^2 \to \ell^2$ defined by
	\begin{equation*} 
		\quad w_n = \sqrt{\frac{n+1}{n}}
	\end{equation*}
	satisfies \eqref{charBwHc} but not \eqref{charBwFhc}, and hence $B_w$ is hypercyclic and non-frequently hypercyclic.
	
	The behaviour of weighted shifts acting on the sequence space $c_0$ is not so straightforward.  Bayart and Grivaux~\cite{BG07} identified weighted shifts $B_w$ acting on $c_0$ that are frequently hypercyclic but not chaotic nor mixing.  However, it was shown by Bonilla and Grosse-Erdmann~\cite{BGE07} that every chaotic weighted shift on $c_0$ is frequently hypercyclic.  
	The following, somewhat technical, characterisation of frequently hypercyclic weighted shifts acting on the space $c_0$ was given in \cite{BR15}.
	\begin{theorem}[Bayart and Ruzsa~\cite{BR15}] \label{thm:fhcBwc0}
		Let $w = (w_n)_{n \in \N_0}$ be a bounded sequence of positive integers.  Then $B_w$ is frequently hypercyclic on $c_0$ if and only if there exists a sequence $(M(p))$ of positive real numbers tending to $+\infty$ and a sequence $(E_p)$ of subsets of $\N_0$ such that
		\begin{enumerate}[label=(\roman*),  itemsep=1ex]
			\item for any $p \geq 1$, $\ldens(E_p) >0$,  \label{thm:fhcBwc0i}
			
			\item for any $p,q \geq 1$, $p \neq q$, $\left( E_p + [0, p] \right) \cap \left( E_p + [0, q] \right) =  \varnothing$,
			
			\item $\displaystyle \lim_{\substack{n \to +\infty \\ n \in E_p + [0,\, p]}} w_1 \cdots w_n = +\infty$,
			
			\item for any $p, q \geq 1$, for any $n \in E_p$ and any $m \in E_q$ with $m > n$, for any $t \in \{ 0, \dotsc, q \}$, 
			\begin{equation*}
				w_1 \cdots w_{m-n+t} \geq M(p)M(q).
			\end{equation*}
		\end{enumerate}
	\end{theorem}
	We note that in \cite{BR15} versions of Theorems \ref{thm:fhcBwlp} and \ref{thm:fhcBwc0} are also proven for the sequence spaces $\ell^p(\Z)$ and $c_0(\Z)$ indexed over the integers.
	
	Recently Charpentier et al.~\cite{CGEM19} extended Theorem \ref{thm:fhcBwlp} to general classes of Fr\'{e}chet sequence spaces. In particular, they identified K\"{o}the sequence spaces such that frequent hypercyclicity of  $B_w$ is equivalent to chaoticity, and examples of K\"{o}the sequence spaces where $B_w$ is frequently hypercyclic but not chaotic.
	
	Examples of frequently hypercyclic operators can also be found among the classes of composition operators acting on function spaces.  
	It turns out that the composition operator $C_\varphi$ acting on the classical Hardy space $H^2(\D)$ is frequently hypercyclic if and only if $C_\varphi$ is hypercyclic \cite{BG06,BG07}. 
	Furthermore, the ever curious family of rank 1 perturbations of unitary operators constructed by Grivaux~\cite{Gri12} are even frequently hypercyclic and chaotic!

	Many results for hypercyclicity have analogues in the frequently hypercyclic case.  For instance, there exists a sufficient condition known as the Frequent Hypercyclicity Criterion.  If an operator $T$ satisfies this criterion then $T$ is frequently hypercyclic, chaotic and mixing (cf.~\cite[Theorem 9.9, Proposition 9.11]{GEP11}).  However, it was shown in \cite{BG07} that there exist frequently hypercyclic weighted shifts on the space $c_0$ that do not satisfy the Frequent Hypercyclicity Criterion.  
	
	The following theorem contrasts further with the hypercyclic case.
	\begin{theorem}[Shkarin~\cite{Shk09}]
		There exist separable infinite-dimensional Banach spaces that do not support frequently hypercyclic operators.
	\end{theorem}
	Analogous to the chaotic case, it was shown in \cite{Shk09} that operators of the form $I + K$, where $K$ is compact or even strictly singular, cannot be frequently hypercyclic. Thus the hereditarily indecomposable Banach spaces constructed in \cite{GM93} and \cite{AH11} do not support frequently hypercyclic operators.

	However, it was shown by de la Rosa et al.~\cite{DLRLGP12} that the chaotic operator from Theorem \ref{thm:RLGP12} is also frequently hypercyclic, which gives that any complex Banach space with an unconditional basis supports a frequently hypercyclic operator.

	Another significant difference to the hypercyclic case is that the set of frequently hypercyclic vectors is of first category.
	\begin{theorem}[Moothathu~\cite{Moo13}, Bayart and Ruzsa~\cite{BR15}, Grivaux and Matheron~\cite{GM14}] \label{thm:fhcVectsMeagre}
		Let $T$ be frequently hypercyclic. Then the set of frequently hypercyclic vectors for $T$ is a meagre set.
	\end{theorem}
	As a consequence of Theorem \ref{thm:fhcVectsMeagre}, the powerful Baire category theorem is not available in the study of frequent hypercyclicity.  Measure-theoretic techniques have emerged as an effective alternative approach and we discuss this topic in Section \ref{sec:ergodicFHC}.
	
	Recently Grivaux~\cite{Gri18} considered frequently hypercyclic vectors with irregular behaviour.  A frequently hypercyclic vector $x \in X$ is said to have an \emph{irregularly visiting orbit} if there exists a nonempty open $U_0 \subset X$ such that the lower density of its return set $\returnSet{x}{U_0}$ is strictly less than its upper density, i.e.
	\begin{equation*}
		\ldens\left( \returnSet{x}{U_0} \right) < \udens\left( \returnSet{x}{U_0} \right),
	\end{equation*}
	or in other words, the return set $\returnSet{x}{U_0}$ has no density.
	It turns out that every operator satisfying a strong form of the Frequent Hypercyclicity Criterion possesses a frequently hypercyclic vector with an irregularly visiting orbit \cite{Gri18}.  Such examples can be found in the article by Grivaux~\cite{Gri17}, who investigated a particular family of hypercyclic operators introduced by Glasner and Weiss~\cite{GW15}. 
	However, it remains an open question whether all frequently hypercyclic operators admit a frequently hypercyclic vector with an irregularly visiting orbit \cite[Question 3.2]{Gri18}.
	
	Finally we mention that two challenging problems in linear dynamics asked whether every chaotic operator is frequently hypercyclic \cite[Question 6.4]{BG07}, and whether the inverse of a frequently hypercyclic operator is still frequently hypercyclic \cite[Question 4.3]{BG06}, \cite[Problem 44]{GMZ16}.  These questions were recently resolved in the negative by Menet~\cite{Men17, Men19b} with the introduction of \emph{operators of C-type}.  We postpone discussion of operators of C-type and these results until Section \ref{sec:CtypeOpers}.

	\subsection{Upper Frequent and Reiterative Hypercyclicity}
	
	Two classes of operators that have recently attracted much interest are the upper frequently and reiteratively hypercyclic operators.  They provide a fine-grained picture of the quantitative differences in the dynamical behaviours that lie between hypercyclicity and frequent hypercyclicity.
	
	The notion of upper, or $\mathcal{U}$-frequent hypercyclicity was introduced by Shkarin~\cite{Shk09}. 
	We say $T \in \Lx$ is $\mathcal{U}$-\emph{frequently hypercyclic} if there exists $x \in X$ such that for any nonempty open subset $U \subset X$, the return set $\nTU$ has positive upper density, i.e.
	\begin{equation*}
		\udens(\nTU)  > 0.
	\end{equation*}
	Such an $x \in X$ is called a $\mathcal{U}$-\emph{frequently hypercyclic vector} for $T$.  
	
	$\mathcal{U}$-frequent hypercyclicity is by definition weaker than frequent hypercyclicity, as illustrated by the example from \cite{BR15} of a weighted shift on the space $c_0$ that is $\mathcal{U}$-frequently hypercyclic but not frequently hypercyclic.  
	
	In further contrast to the frequently hypercyclic case, it was shown in \cite{BR15} that the set of $\mathcal{U}$-frequently hypercyclic vectors for $T$ is comeagre.
	
	\begin{theorem}[Bayart and Ruzsa~\cite{BR15}]  \label{thm:UfhcComeagre}
		Let $T$ be $\mathcal{U}$-frequently hypercyclic.  Then the set of $\mathcal{U}$-frequently hypercyclic vectors for $T$ is comeagre.
	\end{theorem}
	
	However, analogous to the frequently hypercyclic case, the hereditarily indecomposable Banach spaces constructed in \cite{GM93} and \cite{AH11} do not support $\mathcal{U}$-frequently hypercyclic operators \cite{Shk09}.  
	Thus it follows that the set of $\mathcal{U}$-frequently hypercyclic vectors is either empty or comeagre.
	An open question arising from this fact, and stated in \cite{BGE18}, is whether the set of $\mathcal{U}$-frequently hypercyclic vectors (when it exists) is a $G_\delta$-set.
	
	A general criterion that can be used to demonstrate $\mathcal{U}$-frequent hypercyclicity was introduced in \cite{BMPP16}, and Bonilla and Grosse-Erdmann~\cite{BGE18} subsequently introduced a simplification with the following $\mathcal{U}$-Frequent Hypercyclicity Criterion.  
	
	\begin{theorem}[Bonilla and Grosse-Erdmann~\cite{BGE18}] \label{thm:UFHCC}
		Let $T \in \Lx$.  Suppose that there exist dense subsets $X_0, Y_0 \subset X$ and mappings $S_n \colon Y_0 \to X$, $n \geq 0$.  If for any $y \in Y_0$ and $\varepsilon >0$ there exists $A \subset \N$ with $\udens(A) >0$ and $\delta >0$ such that the following hold,
		\begin{enumerate}[label=(\roman*),  itemsep=1ex]
			\item for any $x \in X_0$ there is some $B \subset A$ with $\udens(B) > \delta$ such that for any $n \in B$
			\begin{equation*}
				\norm{T^n x} < \varepsilon,
			\end{equation*}
			
			\item $\displaystyle \sum_{n \in A} S_n y$ converges,
			
			\item for any $m \in A$
			\begin{equation*}
				\norm{T^m \sum_{n \in A} S_n y - y} < \varepsilon,
			\end{equation*}
		\end{enumerate}
		then $T$ is $\mathcal{U}$-frequently hypercyclic.
	\end{theorem}

	It was shown in \cite{BR15} that weighted shifts $B_w$ acting on the spaces $\ell^p$, $1 \leq p < \infty$, are $\mathcal{U}$-frequently hypercyclic if and only if they are frequently hypercyclic, i.e.\ they satisfy \eqref{charBwFhc}.  They also show in \cite{BR15} that replacing the lower density in Theorem \ref{thm:fhcBwc0} \ref{thm:fhcBwc0i} by the condition $\udens(E_p) >0$ gives a characterisation of the $\mathcal{U}$-frequently hypercyclic weighted shifts on $c_0$.

	The following interesting property of invertible frequently hypercyclic operators was identified in \cite{BR15}.
	\begin{theorem}[Bayart and Ruzsa~\cite{BR15}]
		If $T$ is invertible and frequently hypercyclic then its inverse $T^{-1}$ is $\mathcal{U}$-frequently hypercyclic.
	\end{theorem}

	The notion of reiterative hypercyclicity was introduced by B\`{e}s et al.~\cite{BMPP16} when they considered the upper Banach density of  return sets.
	We say $T \in \Lx$ is \emph{reiteratively hypercyclic} if there exists $x \in X$ such that for every nonempty open $U \subset X$, the return set $\nTU$ has positive upper Banach density, i.e.
	\begin{equation*}
		\ubd\left( \nTU \right)  > 0.
	\end{equation*}

	It was proven in \cite{BMPP16} that there exists a reiteratively hypercyclic weighted shift on $c_0$ that is not $\mathcal{U}$-frequently hypercyclic. 
	They also prove that there exists a mixing (and hence hypercyclic) weighted shift on $\ell^p$  that is not reiteratively hypercyclic \cite{BMPP16}.  
	On the other hand, it has been shown that chaos implies reiterative hypercyclicity~\cite{Men17} and reiterative hypercyclicity implies weak mixing \cite{BMPP16}.
	
	\begin{theorem}[B\`{e}s, Menet, Peris and Puig~\cite{BMPP16}]
		Let $X$ be a separable and infinite-dimensional Banach space.  If $T \in \Lx$ is reiteratively hypercyclic then $T$ is weakly mixing.  
	\end{theorem}
	
	\begin{theorem}[Menet~\cite{Men17}]
		Let $X$ be a separable infinite-dimensional Banach space.  If $T \in \Lx$ is chaotic then $T$ is reiteratively hypercyclic. 
	\end{theorem}
	
	If $T$ is reiteratively hypercyclic then it turns out that its set of reiteratively hypercyclic vectors coincides with the set of $T$-hypercyclic vectors \cite{BMPP16}.  
	Some further nice properties of reiterative hypercyclicity are contained in the following theorems, including a positive answer to the  $\mathcal{U}$-frequently and reiteratively hypercyclic analogues of Herrero's problem.
	
	\begin{theorem}[Bonilla and Grosse-Erdmann~\cite{BGE18}] 
		Let $T \in \Lx$ be invertible. If $T$ is reiteratively hypercyclic then so is its inverse.
	\end{theorem}

	\begin{theorem}[Ernst, Esser and Menet~\cite{EEM20}]  \label{thm:invUfhcRhc}
		If $T \in \Lx$ is $\mathcal{U}$-frequently hypercyclic (resp.\ reiteratively hypercyclic), then $T \oplus T$ is $\mathcal{U}$-frequently hypercyclic (resp.\ reiteratively hypercyclic).
	\end{theorem}
	
	It was shown in \cite{BR15} that the weighted shift $B_w$ acting on $\ell^p$ is $\mathcal{U}$-frequently hypercyclic if and only if it is frequently hypercyclic.  This result was subsequently extended in \cite{BMPP16}, where it was shown that every reiteratively hypercyclic weighted shift on $\ell^p(\N)$ and $\ell^p(\Z)$ is frequently hypercyclic.
	It was observed in \cite{BGE18} that replacing upper density by upper Banach density in Theorem \ref{thm:UFHCC} gives a Reiterative Hypercyclicity Criterion.

	
	\begin{figure}[t]
		\centering
		\begin{tikzcd}[row sep=normal, column sep=normal]
			& \textrm{Frequently hypercyclic} \arrow[d,  Rightarrow, blue]  \\
			\textrm{Chaotic} \arrow[dr,  Rightarrow, blue, end anchor=north west]   & \mathcal{U}\textrm{-frequently hypercyclic} 	\ar[d,  Rightarrow,  blue] \\
			\textrm{Mixing} \arrow[dr, Rightarrow, blue, end anchor=north west] &	 \textrm{Reiteratively hypercyclic} \ar[d,  Rightarrow,  blue]  \\
			& \textrm{Weakly mixing} \arrow[d,  Rightarrow, blue]  \\
			&	 \textrm{Hypercyclic} 
		\end{tikzcd}	
		\caption{Relations between various linear dynamical properties.}
		\label{fig:sec4DynProps}
	\end{figure}
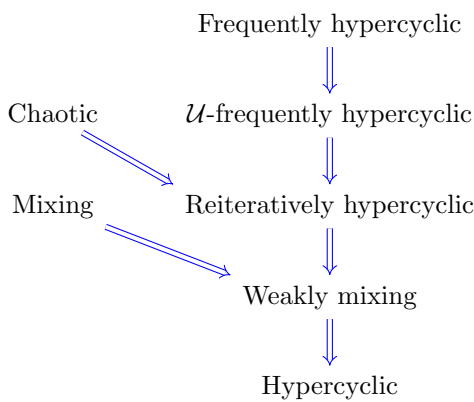
	

	We remark that the studies in \cite{BGE18} and \cite{EEM20} employed the general framework of $\mathcal{A}$-hypercyclicity which was introduced in \cite{BMPP16}, where $\mathcal{A}$ is a Furstenberg family.  
	A nonempty family $\mathcal{A}$ of subsets of $\N_0$ is a \emph{Furstenberg family} if it is hereditary upward, that is if $A \in \mathcal{A}$ and $A \subset B$, then $B \in \mathcal{A}$.  We say $T \in \Lx$ is \emph{$\mathcal{A}$-hypercyclic} if there exists $x \in X$ such that for any nonempty open $U \subset X$ it holds that
	\begin{equation*}
		\nTU \in \mathcal{A}.
	\end{equation*}
	For instance, if $\mathcal{A}$ is the family of nonempty subsets of $\N_0$, then $\mathcal{A}$-hypercyclicity corresponds to hypercyclicity, and if $\mathcal{A}$ is the family of sets with positive lower density then $\mathcal{A}$-hypercyclicity corresponds to frequent hypercyclicity.
	
	With the general machinery of $\mathcal{A}$-hypercyclicity, the results in \cite{BMPP16}, \cite{BGE18} and \cite{EEM20} for $\mathcal{U}$-frequent and reiterative hypercyclicity follow as special cases.
	This approach was also employed by Grosse-Erdmann~\cite{Gro19} to prove the following theorem.
	
	\begin{theorem}[Grosse-Erdmann~\cite{Gro19}]
		Let $B_w$ be an invertible weighted shift on $c_0(\Z)$. If $B_w$ is $\mathcal{U}$-frequently hypercyclic, then so is its inverse $B_w^{-1}$.
	\end{theorem}
	We also note that a study using the framework of Furstenberg families with respect to topological transitivity was conducted by B\`{e}s et al.~\cite{BMPP19}.
	
	Finally we mention that weighted upper and lower densities were studied in Ernst and Mouze~\cite{EM19} and Menet~\cite{Men20} to identify a fine-grained picture of the dynamical notions that lie between $\mathcal{U}$-frequent and frequent hypercyclicity.  In \cite{EEM20} weighted densities were also studied to illustrate the dynamical behaviours lying between reiterative hypercyclicity and $\mathcal{U}$-frequent hypercyclicity.  

	\subsection{Ergodic Theory and Frequent Hypercyclicity}  \label{sec:ergodicFHC}

	In the absence of Baire category techniques, probabilistic methods inspired by ergodic theory have proven a fruitful strategy in the study of frequent hypercyclicity.  This approach is originally due to Flytzanis~\cite{Fly95}, and \cite[Chapter 5]{BM09} contains a comprehensive introduction to this topic.
	
	Unless otherwise stated, in this subsection $X$ denotes a separable complex Banach space with Borel $\sigma$-algebra $\mathcal{B}$.  We let $\mu$ be a finite Borel measure and we let $\T = \{ \lambda \in \C \, : \, \abs{\lambda} = 1\}$ denote the unit circle.  The measure $\mu$ is said to have \emph{full support} if $\mu(A) > 0$ for every nonempty open $A \in \mathcal{B}$.
	
	For a probability space $(X, \mathcal{B}, \mu)$, the measurable map $T \colon (X, \mathcal{B}, \mu) \to (X, \mathcal{B}, \mu)$ is said to be \emph{ergodic} if $\mu$ is $T$-invariant, and if $A \in \mathcal{B}$ satisfies $T^{-1}(A) = A$ up to a set of measure zero then $\mu(A) = 0$ or $1$.  We recall that the measure $\mu$ is said to be $T$-invariant if $\mu\left( T^{-1}(A)\right) = \mu(A)$ for all  $A \in \mathcal{B}$.  
	
	The cornerstone of probabilistic methods in linear dynamics is the following observation.  For $T \in \Lx$, if there exists a $T$-invariant Borel probability measure $\mu$ on $X$, with full support, such that $T$ is an ergodic transformation with respect to $\mu$, then $T$ is frequently hypercyclic (cf.~\cite[Proposition 6.23]{BM09}).
	This can be seen as a consequence of  Birkhoff's ergodic theorem, which states if $T \colon (X, \mathcal{B}, \mu) \to (X, \mathcal{B}, \mu)$ is ergodic with respect to $\mu$, then for any $\mu$-integrable function $f$ on $X$
	\begin{equation}  \label{eq:BirkhoffET}
		\frac{1}{N} \sum_{n=0}^{N-1} f(T^n x) \xrightarrow{N \to \infty} \int_X f \, d\mu 
	\end{equation}
	for $\mu$-almost every $x \in X$ (this is interpreted as meaning that the time average of $f$ with respect to $T$ coincides with its space average).
	The separability of $X$ gives a countable base $\left( V_j \right)$ of open sets of $X$, so when we apply \eqref{eq:BirkhoffET} to the indicator function  $\1_{V_j}$, $j \geq 1$ and using the fact that $\mu$ has full support, we get
	\begin{equation*}
		\frac{1}{N} \sum_{n=0}^{N-1} \1_{V_j} \left(T^n x \right) = \frac{\abs{ \left\lbrace n < N \, : \, T^n x \in V_j \right\rbrace}}{N} \longrightarrow \int_X \1_{V_j} \, d\mu = \mu(V_j) >0.
	\end{equation*}
	Thus there exist subsets $A_j \subset X$, $j \geq 1$, of full measure, such that for any $x \in A_j$
	\begin{equation*}
		\lim_{N \to \infty} \frac{\abs{ \left\lbrace n < N \, : \, T^n x \in V_j \right\rbrace}}{N} \longrightarrow \mu(V_j) >0.
	\end{equation*}
	Since every nonempty open set contains some $V_j$ and since $\bigcap_{j\geq 1} A_j$ has full measure, it follows for $\mu$-almost every $x \in X$ and every nonempty open $U \subset X$ that
	\begin{equation*}
		\liminf_{N \to \infty} \frac{\abs{ \left\lbrace n < N \, : \, T^n x \in U \right\rbrace}}{N}  >0,
	\end{equation*}
	and thus $T$ is frequently hypercyclic.
	
	A practical application of this observation utilises the following fact: if $T$ possesses an abundant supply of unimodular eigenvalues, then $T$ admits an invariant Gaussian measure with respect to which it is ergodic.
	An eigenvector for $T \in \Lx$ is said to be \emph{unimodular} if the associated eigenvalue has modulus 1.  
	We say that $T \in \Lx$ has a \emph{perfectly spanning set of unimodular eigenvectors} if, for every countable $D \subset \T$, we have
	\begin{equation*}
		\cspn{\ker(T - \lambda I) \, : \, \lambda \in \T \setminus D} = X.
	\end{equation*}
	
	Bayart and Grivaux~\cite{BG05} proved for a separable Banach space $X$, that if $T \in \Lx$ has a perfectly spanning set of unimodular eigenvalues then $T$ is hypercyclic.
	They subsequently proved in \cite{BG06} that under these conditions in the Hilbert space setting  $T$ is frequently hypercyclic.  This result was then proven in the Banach space setting by Grivaux~\cite{Gri11}.

	\begin{theorem}[Grivaux~\cite{Gri11}]  \label{thm:perfSpanning}
		For a complex Banach space $X$, if $T \in \Lx$ has a perfectly spanning set of eigenvectors associated to unimodular eigenvalues, then $T$ is frequently hypercyclic. 
	\end{theorem}
	
	Research on this topic initially focused on Gaussian measures, since they are a natural choice in the infinite-dimensional setting. Conditions under which $T \in \Lx$ admits a Gaussian ergodic measure with full support were studied in Bayart and Grivaux~\cite{BG07} and Bayart and Matheron~\cite{BM16}.  However, more recently non-Gaussian measures have been studied by Grivaux~\cite{Gri11}, Murillo-Arcila and Peris~\cite{MP13,MP15}, and Grivaux and Matheron~\cite{GM14}.  We also remark that Theorem \ref{thm:perfSpanning} was proven in \cite{Gri11} for a non-Gaussian measure.   
	
	Some of the main results from \cite{GM14} are the following.  A measure $\mu$ on $X$ is said to be  \emph{continuous} if $\mu(\{ x \}) = 0$ for every $x \in X$.
	
	\begin{theorem}[Grivaux and Matheron~\cite{GM14}]
		Let $X$ be a reflexive Banach space.  Any frequently hypercyclic operator $T \in \Lx$ admits a continuous invariant probability measure with full support.  
	\end{theorem}

	\begin{theorem}[Grivaux and Matheron~\cite{GM14}]
		If $X$ is a reflexive Banach space and if $V \subset X$ is a  nonempty open set, then any frequently hypercyclic operator $T \in \Lx$ admits a continuous ergodic probability measure $\mu$ such that $\mu(V) >0$.
	\end{theorem}
	
	The next theorem concerns the sequence space $c_0(\Z)$ indexed over the integers.
	
	\begin{theorem}[Grivaux and Matheron~\cite{GM14}]  \label{thm:GMnoErgMeas}
		There exists a frequently hypercyclic operator on the space $c_0(\Z)$ that does not admit an ergodic measure with full support.
	\end{theorem}
	
	A key technique employed in \cite{GM14} was the introduction of the parameter $c(T) \in [0, \, 1]$ associated with any hypercyclic $T \in \Lx$.  It is defined for $r>0$ as
	\begin{equation}  \label{defn:cT}
		c(T) = \sup_{x \in \hc{T}} \udens \left( \returnSet{x}{B(0, \, r)}\right),
	\end{equation}
	where $\hc{T}$ is the set of hypercyclic vectors for $T$ and $B(0, \, r)$ is the ball centred at the origin of radius $r>0$.
	It describes the maximal frequency with which the orbit of a hypercyclic vector $x$ visits a ball centred at the origin.
	Its usefulness can be seen in the following theorem.
	\begin{theorem}[Grivaux and Matheron~\cite{GM14}]
		Let $T \in \Lx$.  If $T$ admits an ergodic measure with full support then $c(T) = 1$.
	\end{theorem}
	In \cite{GM14} they show that $\udens \left( \returnSet{x}{B(0, \, r)}\right) = c(T)$ for a comeagre set of hypercyclic vectors and so it follows that
	\begin{equation*}
		c(T) = \sup \left\lbrace c \geq 0 \, : \, \udens \left( \returnSet{x}{B(0, \, r)}\right) \geq c \;\textrm{ for a comeagre set of } x \in \hc{T} \right\rbrace.
	\end{equation*}
	
	Menet~\cite{Men17} recently proved with the following theorem that there exists a chaotic operator that admits only countably many unimodular eigenvalues.  Hence the unimodular eigenvectors cannot be perfectly spanning.  This answered a long-standing question by Flytzanis~\cite{Fly95} that asked whether every hypercyclic operator with unimodular eigenvectors that span a dense subspace necessarily possesses uncountably many unimodular eigenvectors.  
	\begin{theorem}[Menet~\cite{Men17}]  \label{thm:chaoticCtblyUnimod}
		Let $X$ be the complex Banach space $c_0$ or $\ell^p$, for $1 \leq p < \infty$.  
		There exists a chaotic operator $T \in \Lx$ that possesses only countably many unimodular eigenvalues.
	\end{theorem}
	The operator from Theorem \ref{thm:chaoticCtblyUnimod} was an operator of C-type and further discussion of this family of operators appears in Section \ref{sec:CtypeOpers}.

	Another motivation for looking beyond Gaussian measures is that it was recently shown by Grivaux et al.~\cite[Section 2.5]{GMM17} that a \emph{typical} hypercyclic operator $T$ acting on a Hilbert space $H$ does not possess eigenvalues.  
	We will not elaborate here on what is meant in \cite{GMM17} by typical, but we remark that an immediate consequence is that a \emph{typical} hypercyclic operator is not chaotic.

	\section{Li-Yorke and Distributional Chaos}  \label{sec:DistChaos}

	The term \emph{chaos} first appeared in mathematical literature in an article by Li and Yorke~\cite{LY75}, where they studied the dynamical behaviour of interval maps with period three.   Schweizer and Sm{\'\i}tal~\cite{SS94} subsequently introduced the stronger notion of \emph{distributional chaos} for self-maps of a compact interval.  The study of distributional chaos in the linear dynamical setting was initiated by Mart{\'{\i}}nez-Gim{\'e}nez et al.~\cite{MGOP09}.
	
	The operator $T \in \Lx$ is said to be \emph{Li-Yorke chaotic} if there exists an uncountable set $\Gamma \subset X$ such that for each distinct pair $(x, \, y) \in \Gamma \times \Gamma$ we have 
	\begin{equation*}
		\liminf_{n \to \infty} \norm{ \, T^n x - T^n y \,  } = 0 \qquad \textrm{and} \qquad	\limsup_{n \to \infty} \norm{ \,  T^n x - T^n y \, } > 0.
	\end{equation*}
	This definition captures local aspects of the dynamical behaviour of pairs of vectors by describing orbits that are proximal without being asymptotic.
	
	The connection between Li-Yorke chaos and the property of irregularity was identified by Berm{\'u}dez et al.~\cite{BBMGP11}.
	We say that $x \in X$ is an \emph{irregular vector} for $T$  if there exist increasing sequences $(j_k)$ and $(n_k)$ of positive integers such that
	\begin{equation*}
		\lim_{k \to \infty} T^{j_k} x = 0  \quad \textnormal{ and } \quad \lim_{k \to \infty} \norm{ T^{n_k} x } = \infty.
	\end{equation*}
	This notion was introduced by Beauzamy~\cite{Bea88} for Banach spaces and it was generalised to  the Fr\'{e}chet space setting by Bernardes et al.~\cite{BBMP15}.
	
	It was shown in \cite{BBMGP11} that $T$ is Li-Yorke chaotic if and only if $T$ admits an irregular vector.  It is also well known that hypercyclic vectors are irregular.  On the other hand, families of operators that do not admit an irregular vector include compact and normal operators.
	
	Li-Yorke chaotic weighted shifts were characterised in \cite{BBMGP11} with the following theorem.
	\begin{theorem}[Berm{\'u}dez, Bonilla, Mart{\'\i}nez-Gim{\'e}nez and Peris~\cite{BBMGP11}] \label{thm:LiYorkeBwShiftsChar}
		Let $X = c_0$ or $\ell^p$, $1\leq p < \infty$.  The weighted shift $B_w \colon X \to X$ is Li-Yorke chaotic if and only if
		\begin{equation*}
			\sup_{\substack{n \in \N \\ m > n}} \, \prod_{j=n}^m  \abs{w_j} = \infty.
		\end{equation*}
	\end{theorem}
	Comparing Theorem \ref{thm:LiYorkeBwShiftsChar} to \eqref{charBwHc}, we see that it is possible to define a weight sequence such that $B_w$ is Li-Yorke chaotic but not hypercyclic.
	
	The Li-Yorke chaotic composition operators and adjoint multipliers were characterised in \cite{BBMP15}.  In particular, they show for an automorphism $\varphi$ of a domain $\Omega \subset \C$, that the composition operator $C_\varphi$ acting on the space $H(\Omega)$ of holomorphic functions on $\Omega$ is Li-Yorke chaotic if and only if it is hypercyclic.

	A strengthening of Li-Yorke chaos was introduced in \cite{SS94} with the notion of distributional chaos.  
	Before we give the definition of distributional chaos we introduce the following functions.
	Given $\delta > 0$, we define the \emph{lower} and \emph{upper distributional functions} of $x,y \in X$ associated to $T \in \Lx$ as, respectively,
	\begin{align}
		F_{x,y}(\delta) &\coloneqq \ldens\left(\left\lbrace j \in \N \,:\, \norm{ T^j x - T^j y}  < \delta \right\rbrace \right) \label{defn:lowDistFn}
		\shortintertext{and}
		F_{x,y}^*(\delta) &\coloneqq \udens\left(\left\lbrace j \in \N \,:\, \norm{ T^j x - T^j y} < \delta \right\rbrace \right). \label{defn:upDistFn}
	\end{align}
	Note that $F_{x,y}$ and $F_{x,y}^*$ are nondecreasing maps on $(0,\, \infty)$ with $0 \leq F_{x,y} \leq F_{x,y}^* \leq 1$.
	
	If the pair $(x,y)$ satisfy $F^*_{x,y} \equiv 1$ and $F_{x,y}(\varepsilon) = 0$
	for some $\varepsilon > 0$,
	then $(x,y)$ is called a \emph{distributionally chaotic pair}. The operator $T$ is said to be
	\emph{distributionally chaotic}  if there exists an
	uncountable set $\Gamma \subset Y$ such that every distinct pair $(x,y) \in \Gamma \times \Gamma$ is a distributionally chaotic pair for $T$.  The set $\Gamma$ is known as a \emph{distributionally scrambled set} for $T$.
	We say that $T$ is \emph{densely distributionally chaotic} if the set $\Gamma$ may be chosen to be dense in $X$. 
	
	It follows from results in \cite{BBMGP11} and Bernardes et al.~\cite{BBMP13} that in the linear setting  distributional chaos is equivalent to the more tractable notion of distributional irregularity.
	We say $x \in X$ is a \emph{distributionally irregular vector} for $T$
	if there exists $A, B \subset \N$ with 
	\begin{align*}
		\udens(A) = 1 &= \udens(B)
		\shortintertext{such that}
		\lim_{\substack{n \to \infty \\ n \in A}} T^n x = 0  \quad \textnormal{and}& \quad	\lim_{\substack{n \to \infty \\ n \in B}} \norm{ T^n x } = \infty.
	\end{align*}
	This strengthening of irregularity, introduced in \cite{BBMGP11}, describes a more complicated statistical dependence between orbits.


	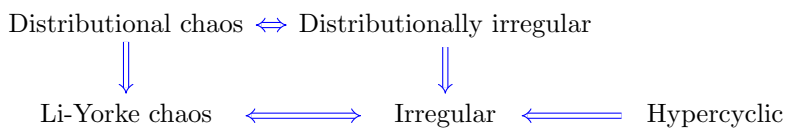
\begin{figure}[t]
		\centering
		\begin{tikzcd}[row sep=normal, column sep=small]
			\textrm{Distributional chaos} 	\ar[d,  Rightarrow,  blue] & \textrm{Distributionally irregular} \arrow[l, Leftrightarrow, blue]  \arrow[d, Rightarrow, blue] & 	\\
			\textrm{Li-Yorke chaos} &  \textrm{Irregular} \arrow[l, Leftrightarrow, blue, shorten >= .8em,shorten <= .8em,] & \textrm{Hypercyclic} \arrow[l, Rightarrow, blue, shorten >= .5em,shorten <= .5em,]	
		\end{tikzcd}	
		\caption{Relations between dynamical properties introduced in Section \ref{sec:DistChaos}.}
		\label{fig:distChaos}
	\end{figure}
	

	Examples of distributionally chaotic operators include the differentiation operator $D$ acting on the space $\entireFns$ of entire functions (this follows from \cite[Corollary 17]{BBMP13}), and weighted  shifts acting on the sequence spaces $\ell^p$ \cite{MGOP09}.  
	
	We have that every infinite-dimensional separable Banach space supports a distributionally chaotic operator.
	\begin{theorem}[Berm\'{u}dez, Bonilla, Mart{\'\i}nez-Gim{\'e}nez and Peris~\cite{BBMGP11}]
		Every infinite-dimensional separable Banach space supports a distributionally chaotic operator.
	\end{theorem}
	
	On the other hand, in \cite[Example 13]{MGOP09} they gave an example of a backward shift acting on a weighted $\ell^p$ space that is distributionally chaotic but not hypercyclic.
	Backward shifts acting on weighted $\ell^p$ spaces are also a source of examples of operators that are mixing but not distributionally chaotic \cite[Theorem 2.1]{MGOP13}.  
	In \cite{MGOP09} they also give an example of a hypercyclic and distributionally chaotic operator that is not chaotic in the sense of Devaney.
	
	In \cite{BR15} they identify, as stated in the theorem below, the existence of frequently hypercyclic operators that are not distributionally chaotic.
	
	\begin{theorem}[Bayart and Ruzsa~\cite{BR15}] 
		There exists a frequently hypercyclic weighted shift acting on $c_0(\Z)$ that is not distributionally chaotic.
	\end{theorem}
	
	However, for unilateral weighted shifts (i.e.\ indexed over $\N$), the properties of frequent hypercyclicity and distributional chaos are equivalent.
	\begin{theorem}[Bayart and Ruzsa~\cite{BR15}]  
		A frequently hypercyclic weighted shift acting on $c_0$ or $\ell^p$, $1 \leq p < \infty$ is distributionally chaotic. 
	\end{theorem}

	For a Banach space $X$, if the set $\{ x \in X \,:\, T^n x \to 0\}$ is dense in $X$, then it follows that Devaney chaos implies distributional chaos, and $\mathcal{U}$-frequent hypercyclicity implies distributional chaos \cite{BBPW18}.
	
	It was also shown in \cite{BBPW18} 
	if there exists a set $S \subset \N$ with $\udens(S) > 0$ such that
	\begin{equation*}
		\sum_{n \in S} \frac{1}{\abs{w_1 \cdots w_n}^p} < \infty 
	\end{equation*}
	then the weighted shift $B_w$ acting on $\ell^p$, $1 \leq p < \infty$, is densely distributionally chaotic.

	We note Bernardes et al.~\cite{BBPW18} conducted an in-depth study of the differing strengths of distributional chaos.  To do this they modified the definition by varying the values of the distributional functions \eqref{defn:lowDistFn} and \eqref{defn:upDistFn}.  Other interesting directions of research include the notion of mean Li-Yorke chaos that has been investigated by Bernardes et al.~\cite{BBP20}.

	The following theorems highlight that the probabilistic techniques described in Section \ref{sec:ergodicFHC} may also shed some light on the nature of distributional chaos.
	
	\begin{theorem}[Bayart and Ruzsa~\cite{BR15}]
		Let $T \in \Lx$ be such that the $\T$-eigenvectors are perfectly spanning with respect to Lebesgue measure. Then $T$ is distributionally chaotic.
	\end{theorem}
	
	It was also shown in \cite[Section 2.7]{GMM17} that typical hypercyclic operator is densely distributionally chaotic.  
	
	\begin{theorem}[Grivaux and Matheron~\cite{GM14}]  
		Let $(X,T)$ be a linear dynamical system, and assume that $T$ admits an ergodic measure with full support. Then $T$ admits a comeagre set of distributionally irregular vectors. This holds in particular if $X$ is a complex Banach space and $T$ has a perfectly spanning set of unimodular eigenvectors.
	\end{theorem}

	\section{Operators of C-Type}  \label{sec:CtypeOpers}
	
	Operators of C-type have emerged as a rich source of chaotic operators that can be fine-tuned to possess very specific properties.  They were introduced by Menet~\cite{Men17} and they provide many important (counter) examples in linear dynamics.  An in-depth study of operators of C-type was conducted by Grivaux et al.~\cite[Section 6]{GMM17} and they have been further developed by Menet~\cite{Men19b, Men19a}.  In this section, unless otherwise specified, \emph{chaos} means chaos in the sense of Devaney.
	
	The class of operators of C-type originated in response to the question, first posed in \cite{BG07}, that asked whether chaotic operators are frequently hypercyclic.  
	We recall in the setting of the sequence spaces $\ell^p$, $1 \leq p < \infty$, that the weighted shift $B_w$ is chaotic if and only if it is frequently hypercyclic \cite{Gro00}, \cite{BR15}.
	On the other hand, the sequence space $c_0$ had hitherto been a fruitful source of counterexamples, indeed it was shown in \cite{BG07} that there exists a frequently hypercyclic weighted shift on $c_0$ that is not chaotic.
	However, in \cite{BGE07} it was proven if $B_w$ is chaotic on $c_0$ then it follows that it is frequently hypercyclic.  There were thus insurmountable obstacles to solving the problem by utilising  established approaches.
	
	The question was finally resolved in \cite{Men17} with the construction of chaotic operators of C-type that are not frequently hypercyclic.   
	An analogous question on whether chaotic operators are distributionally chaotic, posed in \cite{BBMP13},   was also settled in \cite{Men17}.  
	
	\begin{theorem}[Menet~\cite{Men17}]  \label{thm:Men17}
		Let $X$ be the real or complex Banach space $c_0$ or $\ell^p$, for $1 \leq p < \infty$.  
		There exists a chaotic operator $T \in \Lx$ that is not distributionally chaotic nor $\mathcal{U}$-frequently hypercyclic, and hence not frequently hypercyclic.
	\end{theorem}
	
	Before we give the definition of operators of C-type, we introduce some notation.  
	Let $\N_0$ denote the set of non-negative integers and we let $\ell^p(\N_0)$ be the space of $p$-summable sequences indexed over $\N_0$, for $1 \leq p < \infty$.  The canonical basis of $\ell^p(\N_0)$ is denoted by $(e_k)_{k\geq 0}$.  The parameters $\rho, w, \varphi, b$ are defined as follows,
	\begin{itemize}[itemsep=1ex]
		\item $\rho = (\rho_n)_{n\geq 1}$ is a sequence of nonzero complex numbers with $\sum_{n\geq 1} \abs{\rho_n} < \infty$,
		
		\item $w = (w_j)_{j\geq 1}$ is a sequence of complex numbers which is both bounded and bounded below, that is
		\begin{equation*}
			0 < \inf_{k\geq 1} \abs{w_k} \leq \sup_{k \geq 1} \abs{w_k} < \infty,
		\end{equation*}
		
		\item the map $\varphi \colon \N_0 \to \N_0$ is defined such that $\varphi(0) = 0$, $\varphi(n) < n$ for every $n \geq 1$, and the preimage $\varphi^{-1}(l) = \{ n \geq 0 \,:\, \varphi(n) = l \}$ is an infinite set for every $l \geq 0$,   
		
		\item $b = (b_n)_{n \geq 0}$ is a strictly increasing sequence of positive integers such that $b_0 = 0$ and $b_{n+1} - b_n$ is a multiple of $2 \left(b_{\varphi(n) +1} - b_{\varphi(n)} \right)$ for every $n \geq 1$.
	\end{itemize}
	
	\medskip
	
	An \emph{operator of C-type} $T_{\rho,w,\varphi,b} \colon \ell^p(\N_0) \to \ell^p(\N_0)$ is defined as
	\begin{equation*}
		T_{\rho,w,\varphi,b} \, e_k =
		\begin{cases}
			w_{k+1} e_{k+1}, & \quad \textrm{if } k \in [b_n, \, b_{n+1} -1), \, n \geq 0, \\[1ex]
			\rho_n e_{b_{\varphi(n)}} -  \left( \prod\limits_{j=b_n +1}^{b_{n+1} -1} w_j \right)^{-1} e_{b_n}, & \quad  \textrm{if } k = b_{n+1} -1,\, n \geq 1, \\[3ex]
			-  \left( \prod\limits_{j=b_0 +1}^{b_1 -1} w_j \right)^{-1} e_0, &  \quad \textrm{if } k = b_1 -1.
		\end{cases}
	\end{equation*}
	
\begin{figure}[t]
		\centering
		\includegraphics[width=\textwidth]{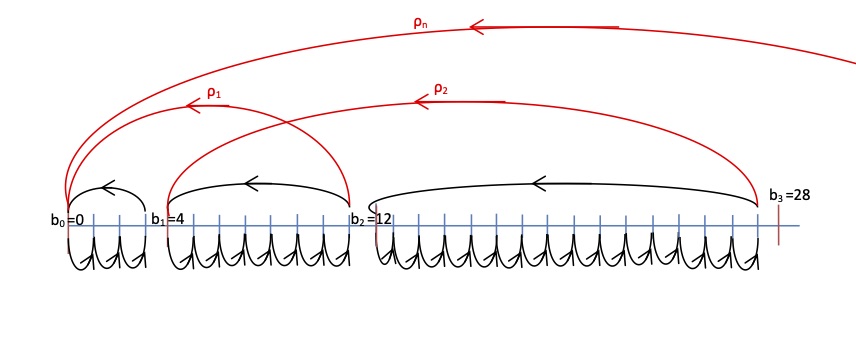}
		\caption{Periodicity of each vector $e_k$. (Figure courtesy of Q.~Menet)} \label{fig:CtypeOper}
	\end{figure}

	A crucial property of operators of C-type is that every basis vector $e_k$ is periodic for $T_{\rho,w,\varphi,b}$.  It is shown in \cite[Fact 6.4]{GMM17} that
	\begin{equation*}
		T_{\rho,w,\varphi,b}^{2(b_{n+1} - b_n)} e_k = e_k, \quad \textrm{if } k \in [b_n, \, b_{n+1}),\, n\geq 0.
	\end{equation*}
	Consequently every finitely supported sequence is periodic for $T_{\rho,w,\varphi,b}$ and thus $T_{\rho,w,\varphi,b}$ possesses a dense set of periodic points.  The periodicity of each $e_k$ is illustrated in Figure \ref{fig:CtypeOper} for a particular sequence $(b_n)$.
	
	It was also shown in \cite[Proposition 6.5]{GMM17} that $T_{\rho,w,\varphi,b}$ is chaotic if
	\begin{equation*}
		\limsup_{\substack{N \to \infty \\ N \in \varphi^{-1}(n)}} \abs{\rho_N} \prod_{j=b_N+1}^{b_{N+1} -1} \abs{w_j} = \infty, \quad \textrm{for every } n \geq 0.
	\end{equation*}
	The term \emph{operator of C-type} was coined in \cite{GMM17}, where the choice of the letter `C' was motivated by their innate connection to \emph{cyclic} and \emph{chaotic} operators.
	
	An inherent characteristic of operators of C-type is an abundant supply of periodic points.  We recall that it is well known for a complex Banach space $X$, that the set of periodic points of $T \in \Lx$ is given by
	\begin{equation*}
		\spn{ x \in X \, : \, Tx = \lambda x, \textrm{ for some root of unity } \lambda \in \C}.
	\end{equation*}
	This is an interesting contrast with the results quoted at the end of Section \ref{sec:ergodicFHC}, which state that a typical hypercyclic operator does not possess any eigenvalues and is not chaotic \cite[Section 2.5]{GMM17}.
	So it is reasonable to ask whether it is possible to find \emph{natural} classes of $\mathcal{U}$-frequently hypercyclic operators  in the Hilbert and Banach space settings that are not frequently hypercyclic.

	
	Another long-standing open problem in linear dynamics, stated in \cite{BG06} and \cite[Problem 44]{GMZ16}, asked whether the inverse of an invertible frequently hypercyclic operator is also frequently hypercyclic.  The answer in the hypercyclic case was already well known, since the inverse of a topologically transitive operator is also topologically transitive (cf.~\cite[Proposition 2.23]{GEP11}).   
	
	As a stepping stone to  resolving this question, Menet~\cite{Men19a} proved the following theorem by constructing a suitable operator of $C$-type.  
	We note that the operators of C-type considered in \cite{GMM17} are not invertible, so a new tweak of the parameters of $T_{\rho,w,\varphi,b}$ was required.
	\begin{theorem}[Menet~\cite{Men19a}]
		There exist invertible $\mathcal{U}$-frequently hypercyclic operators on $\ell^p(\N_0)$, for $1 \leq p < \infty$, such that the inverse is not $\mathcal{U}$-frequently hypercyclic.
	\end{theorem}

	On the other hand, if $T$ is invertible and frequently hypercyclic, it was shown in \cite{BR15} that its inverse $T^{-1}$ is $\mathcal{U}$-frequently hypercyclic.  So a counterexample necessarily requires that $T^{-1}$ is $\mathcal{U}$-frequently hypercyclic but non-frequently hypercyclic.  Such operators were identified in \cite{BR15} and \cite{GMM17}, however they are not invertible.  
	The suitable counterexample was constructed by Menet~\cite{Men19b}.
	
	\begin{theorem}[Menet~\cite{Men19b}]
		There exists an operator on $\ell^1(\N)$ that is invertible and frequently hypercyclic, but with an inverse that is not frequently hypercyclic.
	\end{theorem}

	Another outstanding open question in linear dynamics, originally posed in \cite[Question 4.9]{BG06}, asks whether the direct sum $T \oplus T$ of a frequently hypercyclic operator $T$ is itself frequently hypercyclic.   We recall that there exist hypercyclic operators $T$ such that $T \oplus T$ is not hypercyclic \cite{DLRR09}, \cite{BM07}.  On the other hand, we recall that it was shown in \cite{EEM20} that if $T$ is $\mathcal{U}$-frequently hypercyclic (resp.\ reiteratively hypercyclic), then $T \oplus T$ is  $\mathcal{U}$-frequently hypercyclic (resp.\ reiteratively hypercyclic).
	
	If $T$ is a frequently hypercyclic operator acting on a Banach space $X$, then it is known $T$ is weakly mixing (cf.~\cite[Theorem 9.8]{GEP11}), and thus $T \oplus T$ is hypercyclic on $X \oplus X$.
	This problem is discussed in \cite[Section 6.6]{GMM17}, where they identify a class of operators of C-type that give the following related result.
	\begin{theorem}[Grivaux, Matheron and Menet~\cite{GMM17}]
		Let $p >1$.  There exists a sequence $\left( T_n \right)_{n\geq 1}$ of frequently hypercyclic operators on $\ell^p(\N_0)$ such that the $\ell^p$-sum operator 
		\begin{equation*}
			\bigoplus_{n\geq 1} T_n \colon \bigoplus_{n\geq 1} \ell^p(\N_0) \to \bigoplus_{n\geq 1} \ell^p(\N_0)
		\end{equation*}
		is not $\mathcal{U}$-frequently hypercyclic.
	\end{theorem}

	\subsection{Relations Between Dynamical Properties}

	We briefly revisit Figure \ref{fig:sec4DynProps} to augment it with some of the new examples that emerged following the systematic investigation of operators of C-type conducted in the monster study \cite{GMM17}. 
	
	Most of the implications in Figure \ref{fig:sec4DynProps} follow by definition, apart from the following: chaos implies reiterative hypercyclicity \cite{Men17}, and reiterative hypercyclicity implies weak mixing \cite{BMPP16}.	
	Furthermore, the implications in Figure \ref{fig:sec4DynProps} are strict and in general no further relations hold since there exist examples of operators that are:
	
	\begin{itemize}[itemsep=.7ex]
		\item chaotic but not $\mathcal{U}$-frequently hypercyclic: operators of C-type on $\ell^p$ \cite{Men17} and \cite[Theorem 6.18]{GMM17},
		
		\item chaotic and mixing but not $\mathcal{U}$-frequently hypercyclic: \cite[Example 6.35]{GMM17},
		
		\item chaotic and frequently hypercyclic but not mixing: shift acting on a Hilbert space  \cite{BadGri07},
		
		\item chaotic and frequently hypercyclic but not mixing, nor ergodic with respect to a probability measure with full support: \cite[Example 6.23]{GMM17},
		
		\item frequently hypercyclic but not chaotic nor mixing: weighted shift on $c_0$ \cite{BG07},
		
		\item $\mathcal{U}$-frequently hypercyclic but not frequently hypercyclic: weighted shift on $c_0$ \cite{BR15} and  \cite[Example 6.30]{GMM17},
		
		\item mixing but not reiteratively hypercyclic: weighted shift on $\ell^p$ \cite{BMPP16},
		
		\item reiteratively hypercyclic but not $\mathcal{U}$-frequently hypercyclic: weighted shift on $c_0$ \cite{BMPP16}, 
		
		\item weakly mixing but not mixing: weighted shift on $\ell^p$ \cite{CosSam04},
		
		\item hypercyclic but non-weakly mixing: shift acting on a specific Banach space \cite{DLRR09}, shift acting on $\ell^p$ \cite{BM07}.
	\end{itemize}
	We remark that the results from \cite[Section 6]{GMM17} provide in many instances the first such examples in the Hilbert space setting.

	\section{Hypercyclic Vectors}
	
	A recurring question asked of the author is, `what does a hypercyclic vector look like?'  
	Our primary focus thus far has been on properties of operators, however the vectors that enable complex linear dynamical behaviour also give rise to deep and interesting questions.  In this section we highlight some properties of hypercyclic vectors and we explore the rich algebraic structures contained in the set of hypercyclic vectors.

	We start with the important question of the permissible growth rates of entire functions $f \in \entireFns$ that are hypercyclic with respect to the differentiation operator $D \colon f \mapsto f'$.  
	If $\varphi \colon \R_+ \to \R_+$ is a function with $\varphi(r) \to \infty$ as $r \to \infty$, then there exists a $D$-hypercyclic entire function $f \in \entireFns$ such that
	\begin{equation*}
		\abs{f(z)} \leq \varphi(r) \frac{e^r}{\sqrt{r}}, \quad \textrm{for } \abs{z} = r  \textrm{ sufficiently large}.
	\end{equation*}
	This growth is optimal, since for the critical rate of $e^r/\sqrt{r}$, there does not exist a $D$-hypercyclic entire function $f \in \entireFns$ such that
	\begin{equation*}
		\abs{f(z)} \leq c \frac{e^r}{\sqrt{r}}, \quad \textrm{for } \abs{z} = r >0,
	\end{equation*}
	where $c>0$ is a constant (cf.~\cite[Theorem 4.22]{GEP11}).
	
	As we might expect, the frequently hypercyclic entire functions must grow faster than in the hypercyclic case.  The question of optimal growth in the frequently hypercyclic case was settled by Drasin and Saksman~\cite{DS12} with the careful construction of an entire function with the following growth.
	\begin{theorem}[Drasin and Saksman~\cite{DS12}] \label{thm:DS12}
		For any constant $C >0$, there exists a $D$-frequently hypercyclic entire function $f \in \entireFns$  such that
		\begin{equation*}
			\sup_{\abs{z} = r} \abs{f(z)} \leq C \frac{e^r}{r^{1/4}}, \quad \textrm{ for all } r >0.
		\end{equation*}
		This estimate is optimal, every such function satisfies
		\begin{equation*}
			\limsup_{r \to \infty} r^{1/4} e^{-r} \sup_{\abs{z} = r} \abs{f(z)} >0.
		\end{equation*}
	\end{theorem}
	Analogous questions on the permissible growth rates of $D$-distributionally chaotic entire functions was subsequently investigated in \cite{BGB16} and \cite{GMP18}.

	On the other hand, entire functions $f \in \entireFns$ that are hypercyclic with respect to the translation operators  $T_a \colon f(z) \mapsto f(z +a)$, for $a \neq 0$, may have arbitrarily slow transcendental growth (cf.~\cite[Theorem 4.23]{GEP11}). 
	For a function $\varphi \colon (0,\, \infty) \to [1,\, \infty)$, with $\varphi(r)/r^N \to \infty$ as $r \to \infty$ for any $N \geq 1$, there exists $f \in \entireFns$ that is hypercyclic with respect to  $T_a$ such that
	\begin{equation*}
		\abs{f(z)} \leq C \varphi(r)
	\end{equation*}
	for $\abs{z} = r >0$ and where $C >0$ is a constant.
	
	In the frequently hypercyclic case, Blasco et al.~\cite{BBGE10} identified the following sharp result.
	\begin{theorem}[Blasco, Bonilla and Grosse-Erdmann~\cite{BBGE10}]
		Let the function $\varphi \colon \R_+ \to \R_+$ be arbitrary.  Then there exists a $T_a$-frequently hypercyclic $f \in \entireFns$ with 
		\begin{equation*}
			\sup_{\abs{z} = r} \abs{f(z)} \leq C \varphi(r)
		\end{equation*}
		for $r>0$ sufficiently large and for some constant $C>0$ if and only if 
		\begin{equation*}
			\liminf_{r \to \infty} \frac{\log \varphi(r)}{r} >0.
		\end{equation*}
	\end{theorem}

	Permissible growth rates of hypercyclic, frequently hypercyclic and distributionally chaotic operators acting on the space of harmonic functions on $\R^N$, for $N \geq 2$, have also been studied.  With respect to the partial differentiation operators, the growth of harmonic functions was investigated in \cite{AA98a}, \cite{BBGE10}, \cite{GST19}, \cite{BGB16} and \cite{GMP18}.  
	For translation operators the permissible growth of hypercyclic and frequently hypercyclic harmonic functions was studied in \cite{GCMPR10} and \cite{BBGE10}.

	\subsection{Structural Properties}
	
	We let $\hc{T}$ denote the set of hypercyclic vectors for the hypercyclic operator $T \in \Lx$.  Structural properties of the set $\hc{T}$ have been extensively studied and this has led to some beautiful results.

	For instance, Ansari~\cite{Ans95} proved for a hypercyclic operator $T$ and any $p \in \N$, that the powers $T^p$ are also hypercyclic and that their sets of hypercyclic vectors coincide, i.e.~$\hc{T} = \hc{T^p}$.  This turns out to be a special case of an elegant result by Bourdon and Feldman~\cite{BF03}, who showed that the $T$-orbit of any $x \in X$ is either everywhere dense or nowhere dense in $X$.  
	The analogous result for frequent hypercyclicity was proven in \cite{BG06}, i.e.\ for frequently hypercyclic $T$ and any $p \in \N$, it follows that $T^p$ is frequently hypercyclic and their sets of frequently hypercyclic vectors coincide.
	
	For hypercyclic $T \in \Lx$, it is known that the set $\hc{T}$ is a dense $G_\delta$ subset in $X$ (cf.~\cite[Theorem 1.2]{BM09}), so $\hc{T}$ is large in a topological sense.  As a consequence we have for any hypercyclic $T \in \Lx$, that every vector $x \in X$ can be written as the sum of two hypercyclic vectors, i.e.\ $X = \hc{T} + \hc{T}$ (cf.~\cite[Proposition 1.29]{BM09}).
	
	An important result on the algebraic structure of the set $\hc{T}$, due to Herrero and Bourdon, gives that every hypercyclic operator admits a dense invariant subspace in which every nonzero vector is hypercyclic.  It then follows that the set $\hc{T}$ is a connected subset of $X$ (cf.~\cite[Theorem 2.55, Corollary 2.65]{GEP11}).

	\subsection{Hypercyclic Subspaces}
	
	For a hypercyclic operator $T \in \Lx$, a \emph{hypercyclic subspace} is defined as a closed infinite-dimensional subspace $M \subset X$ such that every nonzero vector in $M$ is hypercyclic for $T$.
	Hypercyclic subspaces have amassed a vast literature since the early work of Bernal Gonz\'{a}lez and Montes-Rodr\'{\i}guez~\cite{BM95}, and introductions to this subject can be found in \cite[Chapter 8]{BM95} and \cite[Chapter 10]{GEP11}.
	
	We briefly remark that this topic fits into the broader area that seeks to identify (inside particular sets) large algebraic structures such as infinite-dimensional vector spaces, closed infinite-dimensional subspaces, and infinitely generated algebras.  These properties are known as lineability, spaceability and algebrability, and we refer the curious reader to \cite{ABPS16} and \cite{BPS14} to learn more about this interesting area.

	In contrast to the generic property that $\hc{T} \cup \{0\}$ always contains a dense subspace of hypercyclic vectors, it turns out that there exist hypercyclic operators that do not admit a hypercyclic subspace. 
	On the other hand,  examples of operators for which every nonzero vector is hypercyclic are uncommon, for instance the construction of Read~\cite{Rea88}.

	Examples of operators that support hypercyclic subspaces include the differentiation and translation operators acting on $\entireFns$ (cf.\ \cite[Examples 10.12 and 10.13]{GEP11}).  For $X = c_0$ or $\ell^p$, $1 \leq p < \infty$, \cite[Example  10.10]{GEP11} recalls that a weighted shift $B_w \in \Lx$ that satisfies the Hypercyclicity Criterion admits a hypercyclic subspace if 
	\begin{equation*}
		\sup_{n\geq 1} \, \limsup_{k \to \infty} \prod_{j=1}^n \abs{w_{j+k}} < \infty.
	\end{equation*}
	However, in \cite[Examples 10.26]{GEP11} it is shown that scalar multiples of the backward shift $cB \in \Lx$, for $\abs{c}>1$, do not possess hypercyclic subspaces.  
	
	The following characterisation of operators that satisfy the Hypercyclicity Criterion and admit a hypercyclic subspace was identified by Gonz{\'a}lez et al.~\cite{GLM00}.
	
	\begin{theorem}[Gonz{\'a}lez, Le{\'o}n-Saavedra and Montes-Rodr{\'{\i}}guez~\cite{GLM00}] \label{thm:HcSubspCharHCC}
		Let $X$ be a separable complex Banach space and suppose that $T \in \Lx$  satisfies the Hypercyclicity Criterion.  The following are equivalent.
		\begin{enumerate}[label=(\roman*),  itemsep=.5ex]
			\item $T$ possesses a hypercyclic subspace.
			
			\item There exists some closed infinite-dimensional subspace $M_0 \subset X$ and an increasing sequence of integers $(n_k)$ such that $T^{n_k}x \to 0$ for all $x \in M_0$.
			
			\item There exists some closed infinite-dimensional subspace $M_0 \subset X$ and an increasing sequence of integers $(n_k)$ such that $\sup_k \lVert T^{n_k}_{\vert M_0} \rVert < \infty$.
			
			\item The essential spectrum of $T$ intersects the closed unit disk.
		\end{enumerate}
	\end{theorem}
	An interesting question arising from Theorem \ref{thm:HcSubspCharHCC} is whether there exists a characterisation of hypercyclic operators that admit a hypercyclic subspace that does not assume the Hypercyclicity Criterion.
	
	We note in the Fr\'{e}chet space setting that the existence of hypercyclic subspaces has recently been investigated by Menet~\cite{Men14, Men15b, Men18}.

	For a family $\left( T_\lambda \right)_{\lambda \in \Lambda}$ of hypercyclic operators acting on the space $X$, we say $x \in X$ is a \emph{common hypercyclic vector} if it is hypercyclic for each $T_\lambda$, $\lambda \in \Lambda$.
	Thus the common hypercyclic vectors are the elements of
	\begin{equation*}
		\bigcap_{\lambda \in \Lambda} \hc{T_\lambda}
	\end{equation*}
	and it immediately follows from the Baire category theorem if $\Lambda$ is countable then the set of common hypercyclic vectors for $\left( T_\lambda \right)_{\lambda \in \Lambda}$ is a dense $G_\delta$-set (cf.~\cite[Proposition 11.2]{GEP11}).

	When the family $\left( T_\lambda \right)_{\lambda \in \Lambda}$ is uncountable, it is recalled in \cite[Chapter 7]{BM09} and \cite[Chapter 11]{GEP11} that common hypercyclic vectors exist for the family of differentiation operators $\left( \lambda D \right)_{\lambda \neq 0}$ for $\lambda \in \C$, and for the family translation operators $\left( T_{e^{i\theta}} \right)_\theta$ for $\theta \in [0,\, 2\pi)$ acting on the space $\entireFns$.  The multiples $\left( \lambda B \right)_\lambda$, for $\abs{\lambda} > 1$, of the backward shift $B$ on the sequence spaces $c_0$ or $\ell^p$, $1\leq p < \infty$, also possess common hypercyclic vectors (cf.~\cite[Example 11.11]{GEP11}). In each of these cases the families of operators in fact admit a dense $G_\delta$ set of common hypercyclic vectors.  
	On the other hand, in \cite[Example 7.1]{BM09} they demonstrate that the hypercyclic weighted shifts $B_w$ on $\ell^2$ do not admit a common hypercyclic vector.

	Hitherto, the results concerning common hypercyclic vectors were primarily for families indexed by a subset of $\R$.  Notable examples in the two-dimensional case were by Shkarin~\cite{Shk10a} and Tsirivas~\cite{Tsi18}.  Bayart~\cite{Bay16} subsequently identified common hypercyclic vectors for higher-dimensional families $\left( T_\lambda \right)_{\lambda \in \R^d}$, for $d \geq 1$.
	
	A natural question arising from this notion asks whether a family of operators can admit a common hypercyclic subspace.  It turns out that even finitely many hypercyclic operators that possess hypercyclic subspaces need not share a common hypercyclic subspace.  For instance, \cite[Example 11.22]{GEP11} demonstrates that the operators $B_w \oplus B_v$ and $B_v \oplus B_w$ acting on the direct sum $\ell^2 \oplus \ell^2$ with
	\begin{equation*}
		w_n = \frac{n+1}{n},\; v_n = 2,\; n \geq 1
	\end{equation*}
	both have hypercyclic subspaces but do not admit a common hypercyclic subspace. 
	
	On the other hand, \cite[Example 11.28]{GEP11} gives for $X = c_0$ or $\ell^p$, $1\leq p < \infty$, that for $T \in \Lx$ defined by 
	\begin{equation*}
		T(x_1,\, x_2,\, x_3, \dotsc) = (x_2,\, x_4,\, x_6,\dotsc),
	\end{equation*}
	the operators $\lambda T$, $\abs{\lambda} > 1$, have a common hypercyclic subspace.
	
	Positive examples also include the countable family of composition operators $C_{\varphi_j}$, for a sequence $(\varphi_j)_{j\geq 1}$ of automorphisms of a domain $\Omega \subset \C$ (cf.\ \cite[p.\ 328]{GEP11}).
	An interesting result from Menet~\cite{Men19} has shown that nonzero multiples of the differentiation operator $D \colon \entireFns \to \entireFns$ and nonzero multiples of translation operators $T_a \colon \entireFns \to \entireFns$, for $a \neq 0$, even share a common hypercyclic subspace!
	
	Further research directions that have been pursued include the study of common $\mathcal{A}$-hypercyclic and common $\mathcal{U}$-frequently hypercyclic vectors, where sufficient conditions for their existence were identified by Mestiri~\cite{Mes20}.  Common and $\mathcal{U}$-frequently hypercyclic subspaces were also studied by B\`{e}s and Menet~\cite{BM15}.
	
	Frequently hypercyclic subspaces are defined as expected. They were originally investigated by Grosse-Erdmann and Bonilla~\cite{BGE12}, who provided the following sufficient condition for their existence.  
	
	\begin{theorem}[Grosse-Erdmann and Bonilla~\cite{BGE12}]
		Let $X$ be a separable Banach space and suppose that $T \in \Lx$ satisfies the Frequent Hypercyclicity Criterion. If there exists a closed infinite-dimensional subspace $M_0 \subset X$ such that for all $x \in X_0$
		\begin{equation*}
			T^n x \to 0
		\end{equation*}
		then $T$ possesses a frequently hypercyclic subspace.
	\end{theorem}
	
	It was demonstrated in \cite{BGE12} that the operator $\varphi(D) \colon \entireFns \to \entireFns$ induced by the  differentiation operator $D$ supports a frequently hypercyclic subspace, where $\varphi$ is an entire function of exponential type that is not a polynomial.  
	
	In contrast, it was proven by Bayart et al.~\cite{BEM16} that the operator $P(D)$ acting on the space $\entireFns$, where $P$ is a nonconstant polynomial, is frequently hypercyclic but does not admit a frequently hypercyclic subspace.  The result from \cite{BEM16} in fact identified a family a frequently hypercyclic operators that admits a hypercyclic subspace but does not support a frequently hypercyclic subspace. 
	They even remark in \cite{BEM16} that it follows from a result of \cite{BM15} that $P(D)$ even admits a $\mathcal{U}$-frequently hypercyclic subspace.  Thus operators of the type $P(D)$ are frequently hypercyclic and support a $\mathcal{U}$-frequently hypercyclic but not a frequently hypercyclic subspace.
	
	The first instance of a frequently hypercyclic operator that supports a hypercyclic but not a frequently hypercyclic subspace was identified by Menet~\cite{Men15} in answer to a question posed in \cite{BGE12}.   
	The example from \cite{Men15} was a weighted shift $B_w \colon \ell^p \to \ell^p$, for $1 \leq p < \infty$, with weights given by $w_k = 2$ if $k \in [a_{2n},\, a_{2n+1})$ for some $n \geq 0$, and $w_k = 1$ otherwise.  Here $(a_n)_{n \geq 0}$ is an appropriately chosen, strictly increasing sequence of integers.

	\subsection{Hypercyclic Algebras}
	
	When the space $X$ supports an algebra structure, a question that has recently attracted much research activity is whether the set of hypercyclic vectors $\hc{T}$ admits, apart from zero, a subalgebra of $X$.  An up-to-date overview of this topic can be found in Bayart et al.~\cite{BCP19}.
	
	We recall that the space $\entireFns$ of entire functions is a Fr\'{e}chet algebra when endowed with pointwise multiplication.  It was shown by Aron et al.~\cite{ACPS07} that non-trivial translation operators $T_a \colon \entireFns \to \entireFns$ do not support a hypercyclic algebras. 
	
	On the other hand, it was independently demonstrated (via different approaches) by Bayart and Matheron~\cite[Section 8.5]{BM09}, and Shkarin~\cite{Shk10b} that the differentiation operator $D \colon \entireFns \to \entireFns$ admits a hypercyclic algebra.  B\`{e}s et al.~\cite{BCP17} subsequently extended this for convolution operators $\Phi(D)$ induced by the differentiation operator.

	\begin{theorem}[B\`{e}s, Conejero and Papathanasiou~\cite{BCP17}]
		Let $\Omega \subset \C$ be a simply connected domain and $H(\Omega)$ the space of holomorphic functions on $\Omega$ endowed with the compact open topology. Let $\Phi$ be a nonconstant polynomial with $\Phi(0) = 0$. Then the set of functions $f \in H(\Omega)$ that generate a hypercyclic algebra for $\Phi(D)$ is comeagre in $H(\Omega)$.
	\end{theorem}
	
	Bayart~\cite{Bay19} identified the following characterisation for convolution operators in the case $\abs{\Phi(0)} < 1$.
	
	\begin{theorem}[Bayart~\cite{Bay19}]
		Let $\Phi$ be a nonconstant entire function of exponential type.
		Assume that $\abs{\Phi(0)} < 1$. The following are equivalent.
		\begin{enumerate}[label=(\roman*),  itemsep=.7ex]
			\item $\Phi(D)$ supports a hypercyclic algebra.
			
			\item  $\Phi$ is not a scalar multiple of an exponential function.
		\end{enumerate}
	\end{theorem}
	
	Sufficient conditions for convolution operators in the case $\abs{\Phi(0)} = 1$ were also given in \cite{Bay19}, and this was extended by B{\`e}s et al.~\cite{BEP20} with the following theorem.
	
	\begin{theorem}[B{\`e}s, Ernst and Prieto~\cite{BEP20}]  \label{thm:phi1hcAlg}
		Let $\Phi$ be a nonconstant entire function of exponential type with $\Phi(0) = 1$.  If $\Phi$ is of subexponential growth, then $\Phi(D)$ has a hypercyclic algebra.
	\end{theorem}
	In \cite{BEP20} they augment Theorem \ref{thm:phi1hcAlg} by giving sufficient conditions for convolution operators to admit a hypercyclic algebra in the case when $\Phi$ is not of subexponential growth.
	
	The picture is not so clear for convolution operators with $\abs{\Phi(0)} > 1$.  An example in \cite{Bay19} demonstrates for $\Phi(z) = 2\exp(-z) + \sin z$, that $\Phi(D)$ supports a hypercyclic algebra.  The following special case was also given in \cite{BCP19}.
	
	\begin{theorem}[Bayart, Costa J\'{u}nior and Papathanasiou~\cite{BCP19}]
		Let $\Phi$ be a nonconstant entire function of exponential type that is not a multiple of an exponential function. Assume that $\abs{\Phi(0)} > 1$ and that there exists some $w \in \C$ such that $\abs{\Phi(tw)} \to 0$ as $t \to +\infty$. Then $\Phi(D)$ supports a hypercyclic algebra.
	\end{theorem}
	
	The question of whether the operator $D$ and convolution operators admit dense and infinitely generated hypercyclic algebras was recently investigated and answered in the affirmative by Bayart~\cite{Bay19}, Bernal-Gonz\'{a}lez and Calder\'{o}n-Moreno~\cite{BC19}, Falc\'{o} and Grosse-Erdmann~\cite{FGE20a}, and B{\`e}s and Papathanasiou~\cite{BP20}.

	For a domain $\Omega \subset \C$, B\`{e}s et al.~\cite{BCP18} proved that weighted composition operators $\wcomp \colon H(\Omega) \to H(\Omega)$ cannot even support a supercyclic algebra, where a supercyclic algebra is defined as expected. 
	However, the following theorem demonstrates that the situation is different if we consider the operator $P(C_\varphi)$ given by the nonconstant polynomial $P$.
	
	\begin{theorem}[Bayart~\cite{Bay19}]
		Let $\Omega \subset \C$ be a simply connected domain and let $\varphi$ be a univalent holomorphic self-map of $\Omega$ which has no fixed point in $\Omega$. Assume $P$ is a nonconstant polynomial that is not a multiple of $z$ and satisfies $\abs{P(1)} < 1$. Then $P(C_\varphi)$ supports a hypercyclic algebra.
	\end{theorem}
	
	The sequence space $\ell^1$ is a Banach algebra when endowed with the convolution product $(u*v)(k) = \sum_{j=0}^k u_j v_{k-j}$, for $u, v \in \ell^1$.
	It was shown in \cite[p.~217]{BM09} that twice the backward shift $2B \colon \ell^1 \to \ell^1$ admits a hypercyclic algebra.  This was subsequently extended in \cite{Bay19} with the following theorem.
	
	\begin{theorem}[Bayart~\cite{Bay19}]
		We consider the backward shift $B \colon \ell^1 \to \ell^1$ and let $P$ be a nonconstant polynomial with complex coefficients.  If $P(\D) \cap \T \neq \varnothing$, then the operator $P(B)$ admits a hypercyclic algebra.
	\end{theorem}
	
	We note that an interesting example from \cite[Example 4.12]{BCP19} identifies an invertible operator $T$ acting on a Banach algebra, such that $T$ supports a hypercyclic algebra while its inverse $T^{-1}$ doe not admit a hypercyclic algebra.
	
	The topic of hypercyclic algebras is rapidly evolving and research is expanding into different directions such as the investigation of frequently hypercyclic algebras \cite{FGE19}, \cite{BCP19}.

	\section{Some Open Questions}
	
	We finish by listing some interesting questions in linear dynamics that remain open.  Most of the questions already appear in the literature and questions \ref{qu:Gri01}-\ref{qu:Gri04} were kindly suggested by Sophie Grivaux.  Further open questions can also be found in \cite{Gri17a} and \cite{Gro17}.

	\begin{qu}  \label{qu:Gri01}
		Find \emph{natural} classes of $\mathcal{U}$-frequently hypercyclic operators on Hilbert and Banach spaces that are not frequently hypercyclic.  Operators of C-type satisfying this description were identified in \cite{GMM17}, but do examples exist among more \emph{classical} operators?
	\end{qu}
	
	\begin{qu}
		Characterise the operators $T_{\lambda, \omega} = D_\lambda + B_\omega$ of the form diagonal plus weighted backward shift, acting on a complex separable infinite-dimensional Hilbert space $\mathcal{H}$, which are hypercyclic, frequently hypercyclic or $\mathcal{U}$-frequently hypercyclic. Here, with respect to a fixed orthonormal basis of $\mathcal{H}$, the diagonal coefficients of $D_\lambda$ are given by $\lambda = (\lambda_k)$, pairwise distinct unimodular complex numbers such that $\lambda_k$ tends to 1 as $k \to \infty$, and for all $M>0$ the weights $\omega = (\omega_k)$ satisfy $0 < \omega_k \leq M$ for every $k \geq 1$. (cf.~\cite[Section 4.3]{GMM17})
	\end{qu}

	\begin{qu}  \label{qu:Gri04}
		On the classical Hardy space $H^2(\D)$, characterise the hypercyclic Toeplitz operators, i.e.\ operators whose matrices with respect to the standard basis of $H^2(\D)$ have constant diagonals (cf.~\cite{BL16}).
	\end{qu}

	\begin{qu}
		Do examples of invertible frequently hypercyclic operators with non-frequently hypercyclic inverse exist on Hilbert or reflexive Banach spaces? (cf.~\cite{Men19b})
	\end{qu}
	
	\begin{qu}
		When the operator $T$ is frequently hypercyclic, does it follow that the direct sum $T \oplus T$ is frequently hypercyclic? (cf.~\cite{BG07})
	\end{qu}
	
	\begin{qu}
		Does every frequently hypercyclic operator admit a frequently hypercyclic vector with an irregularly visiting orbit?  (cf.~\cite[Question 3.2]{Gri18})
	\end{qu}

	\begin{qu}
		Does there exist a characterisation of operators that admit a hypercyclic subspace that does not assume the operator satisfies the Hypercyclicity Criterion?
	\end{qu}

	\section*{Acknowledgements}

	The author is grateful to Fr\'{e}d\'{e}ric Bayart, Sophie Grivaux, Karl Grosse-Erdmann, Quentin Menet and Alfred Peris for helpful suggestions during the initial stages of writing this survey.  He wishes to thank Tom Carroll and Quentin Menet for reading this text and for valuable comments that have improved the article. He is also grateful to the anonymous referee for helpful remarks that improved the text. Permission from Quentin Menet for the use of Figure \ref{fig:CtypeOper} is gratefully acknowledged.  The author would also like to thank the editor Tony O'Farrell for the invitation to write this survey.


\begin{thebibliography}{100}
		
		\bibitem{AA98a}
		M.~P. Aldred and D.~H. Armitage.
		\newblock Harmonic analogues of {G}. {R}. {M}ac{L}ane's universal functions.
		\newblock {\em J. London Math. Soc. (2)}, 57(1):148--156, 1998.
		
		\bibitem{Ans95}
		S.~I. Ansari.
		\newblock Hypercyclic and cyclic vectors.
		\newblock {\em J. Funct. Anal.}, 128(2):374--383, 1995.
		
		\bibitem{Ans97}
		S.~I. Ansari.
		\newblock Existence of hypercyclic operators on topological vector spaces.
		\newblock {\em J. Funct. Anal.}, 148(2):384--390, 1997.
		
		\bibitem{AH11}
		S.~A. Argyros and R.~G. Haydon.
		\newblock A hereditarily indecomposable $\mathcal{L_\infty}$-space that solves
		the scalar-plus-compact problem.
		\newblock {\em Acta Math.}, 206(1):1--54, 2011.
		
		\bibitem{ABPS16}
		R.~M. Aron, L.~Bernal~Gonz\'{a}lez, D.~M. Pellegrino, and J.~B.
		Seoane~Sep\'{u}lveda.
		\newblock {\em Lineability: the search for linearity in mathematics}.
		\newblock Monographs and Research Notes in Mathematics. CRC Press, Boca Raton,
		FL, 2016.
		
		\bibitem{ACPS07}
		R.~M. Aron, J.~A. Conejero, A.~Peris, and J.~B. Seoane-Sep\'{u}lveda.
		\newblock Powers of hypercyclic functions for some classical hypercyclic
		operators.
		\newblock {\em Integr. Equ. Oper. Theory}, 58(4):591--596, 2007.
		
		\bibitem{BadGri07}
		C.~Badea and S.~Grivaux.
		\newblock Unimodular eigenvalues, uniformly distributed sequences and linear
		dynamics.
		\newblock {\em Adv. Math.}, 211(2):766--793, 2007.
		
		\bibitem{BBCDS92}
		J.~Banks, J.~Brooks, G.~Cairns, G.~Davis, and P.~Stacey.
		\newblock On {D}evaney's definition of chaos.
		\newblock {\em Amer. Math. Monthly}, 99(4):332--334, 1992.
		
		\bibitem{BL16}
		A.~Baranov and A.~Lishanskii.
		\newblock Hypercyclic {T}oeplitz operators.
		\newblock {\em Results Math.}, 70(3-4):337--347, 2016.
		
		\bibitem{Bay10}
		F.~Bayart.
		\newblock Parabolic composition operators on the ball.
		\newblock {\em Adv. Math.}, 223(5):1666--1705, 2010.
		
		\bibitem{Bay16}
		F.~Bayart.
		\newblock Common hypercyclic vectors for high-dimensional families of
		operators.
		\newblock {\em Int. Math. Res. Not. IMRN}, (21):6512--6552, 2016.
		
		\bibitem{Bay19}
		F.~Bayart.
		\newblock Hypercyclic algebras.
		\newblock {\em J. Funct. Anal.}, 276(11):3441--3467, 2019.
		
		\bibitem{BC13}
		F.~Bayart and S.~Charpentier.
		\newblock Hyperbolic composition operators on the ball.
		\newblock {\em Trans. Amer. Math. Soc.}, 365(2):911--938, 2013.
		
		\bibitem{BCP19}
		F.~Bayart, F.~Costa~J\'{u}nior, and D.~Papathanasiou.
		\newblock Baire theorem and hypercyclic algebras.
		\newblock arXiv:1910.05000, 2019.
		
		\bibitem{BEM16}
		F.~Bayart, R.~Ernst, and Q.~Menet.
		\newblock Non-existence of frequently hypercyclic subspaces for {$P(D)$}.
		\newblock {\em Israel J. Math.}, 214(1):149--166, 2016.
		
		\bibitem{BG04}
		F.~Bayart and S.~Grivaux.
		\newblock Hypercyclicit\'e: le r\^ole du spectre ponctuel unimodulaire.
		\newblock {\em C. R. Math. Acad. Sci. Paris}, 338(9):703--708, 2004.
		
		\bibitem{BG05}
		F.~Bayart and S.~Grivaux.
		\newblock Hypercyclicity and unimodular point spectrum.
		\newblock {\em J. Funct. Anal.}, 226(2):281--300, 2005.
		
		\bibitem{BG06}
		F.~Bayart and S.~Grivaux.
		\newblock Frequently hypercyclic operators.
		\newblock {\em Trans. Amer. Math. Soc.}, 358(11):5083--5117, 2006.
		
		\bibitem{BG07}
		F.~Bayart and S.~Grivaux.
		\newblock Invariant {G}aussian measures for operators on {B}anach spaces and
		linear dynamics.
		\newblock {\em Proc. Lond. Math. Soc. (3)}, 94(1):181--210, 2007.
		
		\bibitem{BGENP08}
		F.~Bayart, K.-G. Grosse-Erdmann, V.~Nestoridis, and C.~Papadimitropoulos.
		\newblock Abstract theory of universal series and applications.
		\newblock {\em Proc. Lond. Math. Soc. (3)}, 96(2):417--463, 2008.
		
		\bibitem{BM07}
		F.~Bayart and {\'E}.~Matheron.
		\newblock Hypercyclic operators failing the hypercyclicity criterion on
		classical {B}anach spaces.
		\newblock {\em J. Funct. Anal.}, 250(2):426--441, 2007.
		
		\bibitem{BM09}
		F.~Bayart and {\'E}.~Matheron.
		\newblock {\em Dynamics of linear operators}, vol.\ 179 of {\em Cambridge
			Tracts in Mathematics}.
		\newblock Cambridge University Press, Cambridge, 2009.
		
		\bibitem{BM16}
		F.~Bayart and E.~Matheron.
		\newblock Mixing operators and small subsets of the circle.
		\newblock {\em J. Reine Angew. Math.}, 715:75--123, 2016.
		
		\bibitem{BR15}
		F.~Bayart and I.~Z. Ruzsa.
		\newblock Difference sets and frequently hypercyclic weighted shifts.
		\newblock {\em Ergodic Theory Dynam. Systems}, 35(3):691--709, 2015.
		
		\bibitem{Bea88}
		B.~Beauzamy.
		\newblock {\em Introduction to operator theory and invariant subspaces},
		vol.~42 of {\em North-Holland Mathematical Library}.
		\newblock North-Holland Publishing Co., Amsterdam, 1988.
		
		\bibitem{BBMGP11}
		T.~Berm{\'u}dez, A.~Bonilla, F.~Mart{\'\i}nez-Gim{\'e}nez, and A.~Peris.
		\newblock Li-{Y}orke and distributionally chaotic operators.
		\newblock {\em J. Math. Anal. Appl.}, 373(1):83--93, 2011.
		
		\bibitem{Ber99}
		L.~Bernal-Gonz{\'a}lez.
		\newblock On hypercyclic operators on {B}anach spaces.
		\newblock {\em Proc. Amer. Math. Soc.}, 127(4):1003--1010, 1999.
		
		\bibitem{BGB16}
		L.~Bernal-Gonz{\'a}lez and A.~Bonilla.
		\newblock Order of growth of distributionally irregular entire functions for
		the differentiation operator.
		\newblock {\em Complex Var. Elliptic Equ.}, 61(8):1176--1186, 2016.
		
		\bibitem{BC19}
		L.~Bernal-Gonz\'{a}lez and M.~d.~C. Calder\'{o}n-Moreno.
		\newblock Hypercyclic algebras for {$D$}-multiples of convolution operators.
		\newblock {\em Proc. Amer. Math. Soc.}, 147(2):647--653, 2019.
		
		\bibitem{BM95}
		L.~Bernal~Gonz\'{a}lez and A.~Montes-Rodr\'{\i}guez.
		\newblock Universal functions for composition operators.
		\newblock {\em Complex Variables Theory Appl.}, 27(1):47--56, 1995.
		
		\bibitem{BPS14}
		L.~Bernal-Gonz\'{a}lez, D.~Pellegrino, and J.~B. Seoane-Sep\'{u}lveda.
		\newblock Linear subsets of nonlinear sets in topological vector spaces.
		\newblock {\em Bull. Amer. Math. Soc. (N.S.)}, 51(1):71--130, 2014.
		
		\bibitem{BBP20}
		N.~C. Bernardes, A.~Bonilla, and A.~Peris.
		\newblock Mean {L}i-{Y}orke chaos in {B}anach spaces.
		\newblock {\em J. Funct. Anal.}, 278(3):108343, 31, 2020.
		
		\bibitem{BBMP13}
		N.~C. Bernardes, Jr., A.~Bonilla, V.~M{\"u}ller, and A.~Peris.
		\newblock Distributional chaos for linear operators.
		\newblock {\em J. Funct. Anal.}, 265(9):2143--2163, 2013.
		
		\bibitem{BBMP15}
		N.~C. Bernardes, Jr., A.~Bonilla, V.~M{\"u}ller, and A.~Peris.
		\newblock Li-{Y}orke chaos in linear dynamics.
		\newblock {\em Ergodic Theory Dynam. Systems}, 35(6):1723--1745, 2015.
		
		\bibitem{BBPW18}
		N.~C. Bernardes, Jr., A.~Bonilla, A.~Peris, and X.~Wu.
		\newblock Distributional chaos for operators on {B}anach spaces.
		\newblock {\em J. Math. Anal. Appl.}, 459(2):797--821, 2018.
		
		\bibitem{Bes14}
		J.~B{\`e}s.
		\newblock Dynamics of weighted composition operators.
		\newblock {\em Complex Anal. Oper. Theory}, 8(1):159--176, 2014.
		
		\bibitem{BCP17}
		J.~B\`es, J.~A. Conejero, and D.~Papathanasiou.
		\newblock Convolution operators supporting hypercyclic algebras.
		\newblock {\em J. Math. Anal. Appl.}, 445(2):1232--1238, 2017.
		
		\bibitem{BCP18}
		J.~B\`es, J.~A. Conejero, and D.~Papathanasiou.
		\newblock Hypercyclic algebras for convolution and composition operators.
		\newblock {\em J. Funct. Anal.}, 274(10):2884--2905, 2018.
		
		\bibitem{BEP20}
		J.~B\`es, R.~Ernst, and A.~Prieto.
		\newblock Hypercyclic algebras for convolution operators of unimodular constant
		term.
		\newblock {\em J. Math. Anal. Appl.}, 483(1):123595, 25, 2020.
		
		\bibitem{BM15}
		J.~B\`es and Q.~Menet.
		\newblock Existence of common and upper frequently hypercyclic subspaces.
		\newblock {\em J. Math. Anal. Appl.}, 432(1):10--37, 2015.
		
		\bibitem{BMPP16}
		J.~B\`es, Q.~Menet, A.~Peris, and Y.~Puig.
		\newblock Recurrence properties of hypercyclic operators.
		\newblock {\em Math. Ann.}, 366(1-2):545--572, 2016.
		
		\bibitem{BMPP19}
		J.~B\`es, Q.~Menet, A.~Peris, and Y.~Puig.
		\newblock Strong transitivity properties for operators.
		\newblock {\em J. Differential Equations}, 266(2-3):1313--1337, 2019.
		
		\bibitem{BP20}
		J.~B{\`e}s and D.~Papathanasiou.
		\newblock Algebrable sets of hypercyclic vectors for convolution operators.
		\newblock {\em Israel J. Math.}, 238(1):91--119, 2020.
		
		\bibitem{BP99}
		J.~B{\`e}s and A.~Peris.
		\newblock Hereditarily hypercyclic operators.
		\newblock {\em J. Funct. Anal.}, 167(1):94--112, 1999.
		
		\bibitem{Bir20}
		G.~D. Birkhoff.
		\newblock Surface transformations and their dynamical applications.
		\newblock {\em Acta Math.}, 43(1):1--119, 1922.
		
		\bibitem{Bir29}
		G.~D. Birkhoff.
		\newblock D\'emonstration d'un th\'eoreme elementaire sur les fonctions
		entieres.
		\newblock {\em C. R. Acad. Sci. Paris}, 189:473--475, 1929.
		
		\bibitem{BBGE10}
		O.~Blasco, A.~Bonilla, and K.-G. Grosse-Erdmann.
		\newblock Rate of growth of frequently hypercyclic functions.
		\newblock {\em Proc. Edinb. Math. Soc. (2)}, 53(1):39--59, 2010.
		
		\bibitem{BD12}
		J.~Bonet and P.~Doma\'{n}ski.
		\newblock Hypercyclic composition operators on spaces of real analytic
		functions.
		\newblock {\em Math. Proc. Cambridge Philos. Soc.}, 153(3):489--503, 2012.
		
		\bibitem{BMP01}
		J.~Bonet, F.~Mart\'{\i}nez-Gim\'{e}nez, and A.~Peris.
		\newblock A {B}anach space which admits no chaotic operator.
		\newblock {\em Bull. London Math. Soc.}, 33(2):196--198, 2001.
		
		\bibitem{BMP04}
		J.~Bonet, F.~Mart{\'{\i}}nez-Gim{\'e}nez, and A.~Peris.
		\newblock Universal and chaotic multipliers on spaces of operators.
		\newblock {\em J. Math. Anal. Appl.}, 297(2):599--611, 2004.
		
		\bibitem{BP98}
		J.~Bonet and A.~Peris.
		\newblock Hypercyclic operators on non-normable {F}r\'echet spaces.
		\newblock {\em J. Funct. Anal.}, 159(2):587--595, 1998.
		
		\bibitem{BGE07}
		A.~Bonilla and K.-G. Grosse-Erdmann.
		\newblock Frequently hypercyclic operators and vectors.
		\newblock {\em Ergodic Theory Dynam. Systems}, 27(2):383--404, 2007.
		
		\bibitem{BGE12}
		A.~Bonilla and K.-G. Grosse-Erdmann.
		\newblock Frequently hypercyclic subspaces.
		\newblock {\em Monatsh. Math.}, 168(3-4):305--320, 2012.
		
		\bibitem{BGE18}
		A.~Bonilla and K.-G. Grosse-Erdmann.
		\newblock Upper frequent hypercyclicity and related notions.
		\newblock {\em Rev. Mat. Complut.}, 31(3):673--711, 2018.
		
		\bibitem{BF03}
		P.~S. Bourdon and N.~S. Feldman.
		\newblock Somewhere dense orbits are everywhere dense.
		\newblock {\em Indiana Univ. Math. J.}, 52(3):811--819, 2003.
		
		\bibitem{BS90}
		P.~S. Bourdon and J.~H. Shapiro.
		\newblock Cyclic composition operators on {$H^2$}.
		\newblock In {\em Operator theory: operator algebras and applications, {P}art 2
			({D}urham, {NH}, 1988)}, vol.~51 of {\em Proc. Sympos. Pure Math.}, pp.\
		43--53. Amer. Math. Soc., Providence, RI, 1990.
		
		\bibitem{BS97}
		P.~S. Bourdon and J.~H. Shapiro.
		\newblock Cyclic phenomena for composition operators.
		\newblock {\em Mem. Amer. Math. Soc.}, 125(596):x+105, 1997.
		
		\bibitem{CG20}
		T.~Carroll and C.~Gilmore.
		\newblock Weighted composition operators on the {F}ock space and their
		dynamics.
		\newblock arXiv:1911.07254, 2019.
		
		\bibitem{CP11}
		I.~Chalendar and J.~R. Partington.
		\newblock {\em Modern approaches to the invariant-subspace problem}, vol.\ 188
		of {\em Cambridge Tracts in Mathematics}.
		\newblock Cambridge University Press, Cambridge, 2011.
		
		\bibitem{Cha99}
		K.~C. Chan.
		\newblock Hypercyclicity of the operator algebra for a separable {H}ilbert
		space.
		\newblock {\em J. Operator Theory}, 42(2):231--244, 1999.
		
		\bibitem{CGEM19}
		S.~Charpentier, K.~Grosse-Erdmann, and Q.~Menet.
		\newblock Chaos and frequent hypercyclicity for weighted shifts.
		\newblock arXiv:1911.09186, 2019.
		
		\bibitem{CosSam04}
		G.~Costakis and M.~Sambarino.
		\newblock Topologically mixing hypercyclic operators.
		\newblock {\em Proc. Amer. Math. Soc.}, 132(2):385--389, 2004.
		
		\bibitem{DLRLGP12}
		M.~de~la Rosa, L.~Frerick, S.~Grivaux, and A.~Peris.
		\newblock Frequent hypercyclicity, chaos, and unconditional {S}chauder
		decompositions.
		\newblock {\em Israel J. Math.}, 190:389--399, 2012.
		
		\bibitem{DLRR09}
		M.~de~la Rosa and C.~Read.
		\newblock A hypercyclic operator whose direct sum {$T\oplus T$} is not
		hypercyclic.
		\newblock {\em J. Operator Theory}, 61(2):369--380, 2009.
		
		\bibitem{Dev89}
		R.~L. Devaney.
		\newblock {\em An introduction to chaotic dynamical systems}.
		\newblock Addison-Wesley, Redwood City, CA, second edition, 1989.
		
		\bibitem{DS12}
		D.~Drasin and E.~Saksman.
		\newblock Optimal growth of entire functions frequently hypercyclic for the
		differentiation operator.
		\newblock {\em J. Funct. Anal.}, 263(11):3674--3688, 2012.
		
		\bibitem{Enf87}
		P.~Enflo.
		\newblock On the invariant subspace problem for {B}anach spaces.
		\newblock {\em Acta Math.}, 158(3-4):213--313, 1987.
		
		\bibitem{EEM20}
		R.~Ernst, C.~Esser, and Q.~Menet.
		\newblock $\mathcal{U}$-frequent hypercyclicity notions and related weighted
		densities.
		\newblock {\em Israel J. Math.}, to appear, 2019.
		\newblock arXiv:1907.05502.
		
		\bibitem{EM19}
		R.~Ernst and A.~Mouze.
		\newblock A quantitative interpretation of the frequent hypercyclicity
		criterion.
		\newblock {\em Ergodic Theory Dynam. Systems}, 39(4):898--924, 2019.
		
		\bibitem{FGE20a}
		J.~Falc\'{o} and K.-G. Grosse-Erdmann.
		\newblock Algebrability of the set of hypercyclic vectors for backward shift
		operators.
		\newblock {\em Adv. Math.}, 366:107082, 2020.
		
		\bibitem{FGE19}
		J.~Falc\'{o} and K.-G. Grosse-Erdmann.
		\newblock Algebras of frequently hypercyclic vectors.
		\newblock {\em Math. Nachr.}, 293(6):1120--1135, 2020.
		
		\bibitem{FGJ18}
		C.~Fern\'{a}ndez, A.~Galbis, and E.~Jord\'{a}.
		\newblock Dynamics and spectra of composition operators on the {S}chwartz
		space.
		\newblock {\em J. Funct. Anal.}, 274(12):3503--3530, 2018.
		
		\bibitem{Fly95}
		E.~Flytzanis.
		\newblock Unimodular eigenvalues and linear chaos in {H}ilbert spaces.
		\newblock {\em Geom. Funct. Anal.}, 5(1):1--13, 1995.
		
		\bibitem{GS87}
		R.~M. Gethner and J.~H. Shapiro.
		\newblock Universal vectors for operators on spaces of holomorphic functions.
		\newblock {\em Proc. Amer. Math. Soc.}, 100(2):281--288, 1987.
		
		\bibitem{Gil19}
		C.~Gilmore.
		\newblock Dynamics of generalised derivations and elementary operators.
		\newblock {\em Complex Anal. Oper. Theory}, 13(1):257--274, 2019.
		
		\bibitem{GMP18}
		C.~Gilmore, F.~Mart{\'{\i}}nez-Gim{\'e}nez, and A.~Peris.
		\newblock Rate of growth of distributionally chaotic functions.
		\newblock arXiv:1810.09266, 2018.
		
		\bibitem{GST17}
		C.~Gilmore, E.~Saksman, and H.-O. Tylli.
		\newblock Hypercyclicity properties of commutator maps.
		\newblock {\em Integr. Equ. Oper. Theory}, 87(1):139--155, 2017.
		
		\bibitem{GST19}
		C.~Gilmore, E.~Saksman, and H.-O. Tylli.
		\newblock Optimal growth of harmonic functions frequently hypercyclic for the
		partial differentiation operator.
		\newblock {\em Proc.~Roy.~Soc.~Edinburgh Sect.~A}, 149(6):1577--1594, 2019.
		
		\bibitem{GW15}
		E.~Glasner and B.~Weiss.
		\newblock A universal hypercyclic representation.
		\newblock {\em J. Funct. Anal.}, 268(11):3478--3491, 2015.
		
		\bibitem{GS91}
		G.~Godefroy and J.~H. Shapiro.
		\newblock Operators with dense, invariant, cyclic vector manifolds.
		\newblock {\em J. Funct. Anal.}, 98(2):229--269, 1991.
		
		\bibitem{GCMPR10}
		M.~C. G\'omez-Collado, F.~Mart{\'{\i}}nez-Gim{\'e}nez, A.~Peris, and
		F.~Rodenas.
		\newblock Slow growth for universal harmonic functions.
		\newblock {\em J. Inequal. Appl.}, pp.\ 1--6, 2010.
		
		\bibitem{GLM00}
		M.~Gonz{\'a}lez, F.~Le{\'o}n-Saavedra, and A.~Montes-Rodr{\'{\i}}guez.
		\newblock Semi-{F}redholm theory: hypercyclic and supercyclic subspaces.
		\newblock {\em Proc. London Math. Soc. (3)}, 81(1):169--189, 2000.
		
		\bibitem{GM93}
		W.~T. Gowers and B.~Maurey.
		\newblock The unconditional basic sequence problem.
		\newblock {\em J. Amer. Math. Soc.}, 6(4):851--874, 1993.
		
		\bibitem{Gri05}
		S.~Grivaux.
		\newblock Hypercyclic operators, mixing operators, and the bounded steps
		problem.
		\newblock {\em J. Operator Theory}, 54(1):147--168, 2005.
		
		\bibitem{Gri11}
		S.~Grivaux.
		\newblock A new class of frequently hypercyclic operators.
		\newblock {\em Indiana Univ. Math. J.}, 60(4):1177--1201, 2011.
		
		\bibitem{Gri12}
		S.~Grivaux.
		\newblock A hypercyclic rank one perturbation of a unitary operator.
		\newblock {\em Math. Nachr.}, 285(5-6):533--544, 2012.
		
		\bibitem{Gri17}
		S.~Grivaux.
		\newblock Some new examples of universal hypercyclic operators in the sense of
		{G}lasner and {W}eiss.
		\newblock {\em Trans. Amer. Math. Soc.}, 369(11):7589--7629, 2017.
		
		\bibitem{Gri17a}
		S.~Grivaux.
		\newblock Ten questions in linear dynamics.
		\newblock In {\em \'{E}tudes op\'{e}ratorielles}, vol.\ 112 of {\em Banach
			Center Publ.}, pp.\ 143--151. Polish Acad. Sci. Inst. Math., Warsaw, 2017.
		
		\bibitem{Gri18}
		S.~Grivaux.
		\newblock Frequently hypercyclic operators with irregularly visiting orbits.
		\newblock {\em J. Math. Anal. Appl.}, 462(1):542--553, 2018.
		
		\bibitem{GM14}
		S.~Grivaux and {\'E}.~Matheron.
		\newblock Invariant measures for frequently hypercyclic operators.
		\newblock {\em Adv. Math.}, 265:371--427, 2014.
		
		\bibitem{GMM17}
		S.~{Grivaux}, {\'E}.~{Matheron}, and Q.~{Menet}.
		\newblock {Linear dynamical systems on Hilbert spaces:~typical properties and
			explicit examples}.
		\newblock {\em Mem. Amer. Math. Soc.}, to appear, 2017.
		\newblock arXiv:1703.01854.
		
		\bibitem{GR08}
		S.~Grivaux and M.~Roginskaya.
		\newblock On {R}ead's type operators on {H}ilbert spaces.
		\newblock {\em Int. Math. Res. Not. IMRN} Art. ID rnn 083, 42, 2008.
		
		\bibitem{GR14}
		S.~Grivaux and M.~Roginskaya.
		\newblock A general approach to {R}ead's type constructions of operators
		without non-trivial invariant closed subspaces.
		\newblock {\em Proc. Lond. Math. Soc. (3)}, 109(3):596--652, 2014.
		
		\bibitem{Gro99}
		K.-G. Grosse-Erdmann.
		\newblock Universal families and hypercyclic operators.
		\newblock {\em Bull. Amer. Math. Soc. (N.S.)}, 36(3):345--381, 1999.
		
		\bibitem{Gro00}
		K.-G. Grosse-Erdmann.
		\newblock Hypercyclic and chaotic weighted shifts.
		\newblock {\em Studia Math.}, 139(1):47--68, 2000.
		
		\bibitem{Gro17}
		K.-G. Grosse-Erdmann.
		\newblock Frequently hypercyclic operators: recent advances and open problems.
		\newblock In {\em Advanced courses of mathematical analysis {VI}}, pp.\
		173--190. World Sci. Publ., Hackensack, NJ, 2017.
		
		\bibitem{Gro19}
		K.-G. Grosse-Erdmann.
		\newblock Frequently hypercyclic bilateral shifts.
		\newblock {\em Glasg. Math. J.}, 61(2):271--286, 2019.
		
		\bibitem{GEP11}
		K.-G. Grosse-Erdmann and A.~Peris~Manguillot.
		\newblock {\em Linear chaos}.
		\newblock Universitext. Springer, London, 2011.
		
		\bibitem{GMZ16}
		A.~J. Guirao, V.~Montesinos, and V.~Zizler.
		\newblock {\em Open problems in the geometry and analysis of {B}anach spaces}.
		\newblock Springer, 2016.
		
		\bibitem{Her91}
		D.~A. Herrero.
		\newblock Limits of hypercyclic and supercyclic operators.
		\newblock {\em J. Funct. Anal.}, 99(1):179--190, 1991.
		
		\bibitem{Kit82}
		C.~Kitai.
		\newblock {\em Invariant closed sets for linear operators}.
		\newblock PhD thesis, University of Toronto, 1982.
		
		\bibitem{LY75}
		T.~Y. Li and J.~A. Yorke.
		\newblock Period three implies chaos.
		\newblock {\em Amer. Math. Monthly}, 82(10):985--992, 1975.
		
		\bibitem{Mac52}
		G.~R. MacLane.
		\newblock Sequences of derivatives and normal families.
		\newblock {\em J. Analyse Math.}, 2:72--87, 1952.
		
		\bibitem{MGOP09}
		F.~Mart{\'\i}nez-Gim{\'e}nez, P.~Oprocha, and A.~Peris.
		\newblock Distributional chaos for backward shifts.
		\newblock {\em J. Math. Anal. Appl.}, 351(2):607--615, 2009.
		
		\bibitem{MGOP13}
		F.~Mart{\'\i}nez-Gim{\'e}nez, P.~Oprocha, and A.~Peris.
		\newblock Distributional chaos for operators with full scrambled sets.
		\newblock {\em Math. Z.}, 274(1-2):603--612, 2013.
		
		\bibitem{Men14}
		Q.~Menet.
		\newblock Hypercyclic subspaces and weighted shifts.
		\newblock {\em Adv. Math.}, 255:305--337, 2014.
		
		\bibitem{Men15}
		Q.~Menet.
		\newblock Existence and non-existence of frequently hypercyclic subspaces for
		weighted shifts.
		\newblock {\em Proc. Amer. Math. Soc.}, 143(6):2469--2477, 2015.
		
		\bibitem{Men15b}
		Q.~Menet.
		\newblock Hereditarily hypercyclic subspaces.
		\newblock {\em J. Operator Theory}, 73(2):385--405, 2015.
		
		\bibitem{Men17}
		Q.~Menet.
		\newblock Linear chaos and frequent hypercyclicity.
		\newblock {\em Trans. Amer. Math. Soc.}, 369(7):4977--4994, 2017.
		
		\bibitem{Men18}
		Q.~Menet.
		\newblock Invariant subspaces for non-normable {F}r\'{e}chet spaces.
		\newblock {\em Adv. Math.}, 339:495--539, 2018.
		
		\bibitem{Men19}
		Q.~Menet.
		\newblock Existence of common hypercyclic subspaces for the derivative operator
		and the translation operators.
		\newblock {\em Rev. R. Acad. Cienc. Exactas F\'{\i}s. Nat. Ser. A Mat. RACSAM},
		113(2):487--505, 2019.
		
		\bibitem{Men19b}
		Q.~Menet.
		\newblock Inverse of frequently hypercyclic operators.
		\newblock arXiv:1910.04452, 2019.
		
		\bibitem{Men20}
		Q.~Menet.
		\newblock A bridge between $\mathcal{U}$-frequent hypercyclicity and frequent
		hypercyclicity.
		\newblock {\em J. Math. Anal. Appl.}, 482(2):123569, 15, 2020.
		
		\bibitem{Men19a}
		Q.~Menet.
		\newblock Inverse of {$\mathcal{U}$}-frequently hypercyclic operators.
		\newblock {\em J. Funct. Anal.}, 279(4):108543, 2020.
		
		\bibitem{Mes20}
		M.~Mestiri.
		\newblock Common upper frequent hypercyclicity.
		\newblock {\em Studia Math.}, 250(1):1--18, 2020.
		
		\bibitem{Moo13}
		T.~K.~S. Moothathu.
		\newblock Two remarks on frequent hypercyclicity.
		\newblock {\em J. Math. Anal. Appl.}, 408(2):843--845, 2013.
		
		\bibitem{MP13}
		M.~Murillo-Arcila and A.~Peris.
		\newblock Strong mixing measures for linear operators and frequent
		hypercyclicity.
		\newblock {\em J. Math. Anal. Appl.}, 398(2):462--465, 2013.
		
		\bibitem{MP15}
		M.~Murillo-Arcila and A.~Peris.
		\newblock Strong mixing measures for {$C_0$}-semigroups.
		\newblock {\em Rev. R. Acad. Cienc. Exactas F\'{\i}s. Nat. Ser. A Mat. RACSAM},
		109(1):101--115, 2015.
		
		\bibitem{Pal14}
		J.~Pal.
		\newblock Zwei kleine bemerkungen.
		\newblock {\em Tohoku Math. J.}, 6(6):42--–43, 1914.
		
		\bibitem{RR03}
		H.~Radjavi and P.~Rosenthal.
		\newblock {\em Invariant subspaces}.
		\newblock Dover Publications, Inc., Mineola, NY, second edition, 2003.
		
		\bibitem{Rea85}
		C.~J. Read.
		\newblock A solution to the invariant subspace problem on the space {$l_1$}.
		\newblock {\em Bull. London Math. Soc.}, 17(4):305--317, 1985.
		
		\bibitem{Rea88}
		C.~J. Read.
		\newblock The invariant subspace problem for a class of {B}anach spaces. {II}.
		{H}ypercyclic operators.
		\newblock {\em Israel J. Math.}, 63(1):1--40, 1988.
		
		\bibitem{Rez11}
		H.~Rezaei.
		\newblock Chaotic property of weighted composition operators.
		\newblock {\em Bull. Korean Math. Soc.}, 48(6):1119--1124, 2011.
		
		\bibitem{Rol69}
		S.~Rolewicz.
		\newblock On orbits of elements.
		\newblock {\em Studia Math.}, 32:17--22, 1969.
		
		\bibitem{Sal95}
		H.~N. Salas.
		\newblock Hypercyclic weighted shifts.
		\newblock {\em Trans. Amer. Math. Soc.}, 347(3):993--1004, 1995.
		
		\bibitem{SS94}
		B.~Schweizer and J.~Sm{\'\i}tal.
		\newblock Measures of chaos and a spectral decomposition of dynamical systems
		on the interval.
		\newblock {\em Trans. Amer. Math. Soc.}, 344(2):737--754, 1994.
		
		\bibitem{Shk09}
		S.~Shkarin.
		\newblock On the spectrum of frequently hypercyclic operators.
		\newblock {\em Proc. Amer. Math. Soc.}, 137(1):123--134, 2009.
		
		\bibitem{Shk10c}
		S.~Shkarin.
		\newblock A hypercyclic finite rank perturbation of a unitary operator.
		\newblock {\em Math. Ann.}, 348(2):379--393, 2010.
		
		\bibitem{Shk10b}
		S.~Shkarin.
		\newblock On the set of hypercyclic vectors for the differentiation operator.
		\newblock {\em Israel J. Math.}, 180:271--283, 2010.
		
		\bibitem{Shk10a}
		S.~Shkarin.
		\newblock Remarks on common hypercyclic vectors.
		\newblock {\em J. Funct. Anal.}, 258(1):132--160, 2010.
		
		\bibitem{Tan04}
		M.~Taniguchi.
		\newblock Chaotic composition operators on the classical holomorphic spaces.
		\newblock {\em Complex Var. Theory Appl.}, 49(7-9):529--538, 2004.
		
		\bibitem{Tsi18}
		N.~Tsirivas.
		\newblock Existence of common hypercyclic vectors for translation operators.
		\newblock {\em J. Operator Theory}, 80(2):257--294, 2018.
		
	\end{thebibliography}

%
%
	
\end{document}